\documentstyle{amltd}
\begin{document}
\annalsline{151}{2000}
\received{March 12, 1998}
\startingpage{151}
\def\eqref#1{(\ref{#1})}
\def\bye{\end{document}}
 \font\tenrm=cmr10
\input amssym.def
\input amssym.tex
\catcode`\@=11
\font\twelvemsb=msbm10 scaled 1100
\font\tenmsb=msbm10
\font\ninemsb=msbm10 scaled 800
\newfam\msbfam
\textfont\msbfam=\twelvemsb  \scriptfont\msbfam=\ninemsb
  \scriptscriptfont\msbfam=\ninemsb
\def\msb@{\hexnumber@\msbfam}
\def\Bbb{\relax\ifmmode\let\next\Bbb@\else
 \def\next{\errmessage{Use \string\Bbb\space only in math
mode}}\fi\next}
\def\Bbb@#1{{\Bbb@@{#1}}}
\def\Bbb@@#1{\fam\msbfam#1}
\catcode`\@=12

 \catcode`\@=11
\font\twelveeuf=eufm10 scaled 1100
\font\teneuf=eufm10
\font\nineeuf=eufm7 scaled 1100
\newfam\euffam
\textfont\euffam=\twelveeuf  \scriptfont\euffam=\teneuf
  \scriptscriptfont\euffam=\nineeuf
\def\euf@{\hexnumber@\euffam}
\def\frak{\relax\ifmmode\let\next\frak@\else
 \def\next{\errmessage{Use \string\frak\space only in math
mode}}\fi\next}
\def\frak@#1{{\frak@@{#1}}}
\def\frak@@#1{\fam\euffam#1}
\catcode`\@=12
\newcommand{\MA}{Monge-Amp\`ere }
\renewcommand{\c}{\overline}
\newcommand{\cz}{\c{z}}
\newcommand{\pa}{\partial}
\newcommand{\wt}{\widetilde}
\newcommand{\wh}{\widehat}
\newcommand{\e}{\varepsilon}

\newcommand{\qtext}[1]{\quad\hbox{#1}\ \ }
\newcommand{\R}[2]{R^{(#1,#2)}}
\newcommand{\cj}{{\c\jmath}}
\newcommand{\C}{{\Bbb C}}

\newcommand{\calA}{{\cal A}}
\newcommand{\calC}{{\cal C}}
\newcommand{\calE}{{\cal E}}
\newcommand{\calF}{{\cal F}}
\newcommand{\calH}{{\cal H}}
\newcommand{\calJ}{{\cal J}}
\newcommand{\calL}{{\cal L}}
\newcommand{\calM}{{\cal M}}
\newcommand{\calN}{{\cal N}}
\newcommand{\calP}{{\cal P}}
\newcommand{\calQ}{{\cal Q}}
\newcommand{\calR}{{\cal R}}
\newcommand{\calS}{{\cal S}}
\newcommand{\calT}{{\cal T}}

\newcommand{\bA}{{\bf A}}
\newcommand{\bC}{{\bf C}}
\newcommand{\bS}{{\bf S}}
\newcommand{\bT}{{\bf T}}

\newcommand{\su}{{\frak s\frak u}}

\renewcommand{\Re}{{\rm Re}}
\renewcommand{\Im}{{\rm Im}}
\renewcommand{\O}{{\rm O}}
\newcommand{\Ric}{{\rm Ricci}}
\newcommand{\contr}{{\rm contr}}
\newcommand{\tr}{{\rm tr}}
\newcommand{\sig}{{\rm sig}}
\newcommand{\id}{{\rm id}}
\newcommand{\SU}{{\rm SU}}
\newcommand{\U}{{\rm U}}
\newcommand{\F}{{{\rm F}}}
\newcommand{\G}{{{\rm G}}}
\newcommand{\Range}{{\rm Range}}

\newcommand{\contraction}{%
    \hskip2pt\vbox{\hbox{\phantom{2}{\vrule}}\hrule}\hskip2pt}


\title{Construction of boundary invariants\\ and the 
 logarithmic singularity\\ of the Bergman kernel} 
\shorttitle{Boundary invariants and the Bergman kernel} 
 \acknowledgements{This research was supported by Grant-in-Aid  for Scientific Research, 
The Ministry of Education,  Science and Culture, Japan and by NSF grant
\#DMS-9022140 at MSRI.}

\author{Kengo Hirachi}
 
 \institutions{
Osaka University,
Toyonaka, Osaka 560, Japan\\
{\eightpoint {\it E-mail address\/}: hirachi@math.sci.osaka-u.ac.jp}}

\intro
This paper studies Fefferman's program \cite{F3} of expressing the
singularity of the Bergman kernel, for smoothly bounded strictly
pseudoconvex domains $\Omega\subset\C^n$, in terms of local biholomorphic
invariants of the boundary. By \cite{F1}, the Bergman kernel on the
diagonal $K(z,\c z)$ is written in the form
$$
 K=\varphi\, r^{-n-1}+\psi \log r \qtext{with}
  \varphi,\psi\in C^\infty(\c\Omega),
$$ 
where $r$ is a (smooth) defining function of $\Omega$. Recently, Bailey,
Eastwood and Graham \cite{BEG}, building on Fefferman's earlier work
\cite{F3}, obtained a full invariant expression of the strong singularity
$\varphi\, r^{-n-1}$. The purpose of this paper is to give a full invariant
expression of the weak singularity $\psi\log r$.

Fefferman's program is modeled on the heat kernel asymptotics for
Riemannian manifolds,
$$
 K_t(x,x)\sim t^{-n/2}\sum_{j=0}^\infty a_j(x)\,t^j
 \quad \hbox{as} \ t \to+0,
$$ 
in which case the coefficients $a_j$ are expressed, by the Weyl
invariant theory, in terms of the Riemannian curvature tensor and its
covariant derivatives. The Bergman kernel's counterpart of the time
variable $t$ is a defining function $r$ of the domain $\Omega$. By
\cite{F1} and \cite{BS}, the formal singularity of $K$ at a boundary point
$p$ is uniquely determined by the Taylor expansion of $r$ at $p$. Thus one
has hope of expressing $\varphi$ modulo $O^{n+1}(r)$ and $\psi$ modulo
$O^{\infty}(r)$ in terms of local biholomorphic invariants of the boundary,
provided $r$ is appropriately chosen. In \cite{F3}, Fefferman proposed to
find such expressions by reducing the problem to an algebraic one in
invariant theory associated with CR geometry, and indeed expressed
$\varphi$ modulo $O^{n-19}(r)$ invariantly by solving the reduced problem
partially. The solution in \cite{F3} was then completed in \cite{BEG} to
give a full invariant expression of $\varphi$ modulo $O^{n+1}(r)$, but the
reduction is still obstructed at finite order so that the procedure does
not apply to the log term $\psi$. We thus modify the invariant-theoretic
problem in \cite{F3}, \cite{BEG} and solve the modified problem to extend
the reduction.

In the heat kernel case, the reduction to the algebraic problem is done by
using normal coordinates, and the coefficient functions $a_k$ at a point of
reference are $\O(n)$-invariant polynomials in jets of the metric. The CR
geometry counterpart of the normal coordinates has been given by Moser
\cite{CM}. If $\pa\Omega\in C^\omega$ (real-analytic) then, after a change
of local coordinates, $\pa\Omega$ is locally placed in Moser's normal form:
\begin{equation}\label{mnf-rho}
 N(A): \ \
 \rho(z,\c z)=2u-|z'|^2-\sum_{|\alpha|,|\beta|\ge2,l\ge0}
  A_{\alpha\c\beta}^l z'{}^\alpha\c z'{}^\beta v^l=0,
\end{equation} 
where $z'=(z^1,\dots,z^{n-1})$, $z^n=u+i v$, $A=(A_{\alpha\c\beta}^l)$, and
the coefficients $A_{\alpha\c\beta}^l$ satisfy trace conditions which are
linear (see Section 3). For each $p\in\pa\Omega$, Moser's local coordinate
system as above is uniquely determined up to an action of a parabolic
subgroup $H$ of $\SU(1,n)$. Thus $H$-invariant functions of $A$ give rise to
local biholomorphic invariants at the point $p$. Among these invariants, we
define {\it CR invariants of weight $w$} to be polynomials $I(A)$ in $A$
such that
\begin{equation}\label{CRinv-transf}
  I(\wt A)=|\det\Phi'(0)|^{-2w/(n+1)}I(A)
\end{equation} 
for biholomorphic maps $\Phi$ such that $\Phi(0)=0$ and 
$\Phi(N(A))=N(\wt A)$. A CR invariant $I(A)$ defines an assignment, to each
strictly pseudoconvex hypersurface $M\in C^\omega$, of a function 
$I_M\in C^\omega(M)$, which is also called a CR invariant. Here $I_M(p)$,
$p\in M$, is given by taking  a biholomorphic map such that $\Phi(p)=0$,
$\Phi(M)=N(A)$ and then setting
\begin{equation}\label{CRinv-functional}
 I_M(p)=|\det\Phi'(p)|^{2w/(n+1)}I(A).
\end{equation} 
This value is independent of the choice of $\Phi$ with $N(A)$ because of
\eqref{CRinv-transf}. If $M\in C^\infty$ then \eqref{CRinv-functional} gives
$I_M\in C^\infty(M)$, though a normal form of $M$ can be a formal surface.

The difficulty of the whole problem comes from the ambiguity of the choice
of defining functions $r$, and this has already appeared in the problem for
$\varphi$, that is, the problem of finding an expression for $\varphi$ of
the form
\begin{equation}\label{exp-varphi}
 \varphi=\sum_{j=0}^n \varphi_j\,r^j+O^{n+1}(r)
 \quad \hbox{with} \ \
 \varphi_j\in C^\infty(\overline\Omega),
\end{equation} 
such that the boundary value $\varphi_j|_{\pa\Omega}$ is a CR invariant of
weight $j$. Though this expansion looks similar to that of the heat kernel,
the situation is much more intricate. It is impossible to choose an exactly
invariant defining function $r$, and thus the extension of CR invariants
$\varphi_j|_{\pa\Omega}$ to $\Omega$ near $\pa\Omega$, which is crucial, is
inevitably approximate. Fefferman \cite{F3} employed an approximately
invariant defining function $r=r^{\F}$, which was constructed in \cite{F2}
as a smooth approximate solution to the (complex) \MA equation (with zero
Dirichlet condition). This defining function is uniquely determined with
error of order $n+2$ along the boundary, and approximately invariant under
biholomorphic maps $\Phi\colon\Omega\to\wt\Omega$ in the sense that
\begin{equation}\label{trans-r0}
  \wt r\circ\Phi=|\det\Phi'|^{2/(n+1)}r+O^{n+2}(r),
\end{equation} 
for $r=r^{\F}$ and $\wt r=\wt r^{\F}$ associated with $\Omega$ and
$\wt\Omega$, respectively. The defining function $r=r^{\F}$ was used by
\cite{F3} and \cite{BEG} also in the ambient metric construction of the
coefficient functions $\varphi_j$ explained as follows. Let $g[r]$ be the
Lorentz-K\"ahler metric on $\C^*\times\c\Omega\subset\C^{n+1}$ near
$\C^*\times\pa\Omega$ defined by the potential $|z^0|^2r$ ($z^0\in\C^*$).
Then scalar functions are obtained as complete contractions of tensor
products of covariant derivatives of the curvature tensor of $g[r]$. By
\cite{F3} and \cite{BEG}, such complete contractions generate all CR
invariants of weight $\le n$, and each $\varphi_j$ in the expansion of
$\varphi$ is realized by linear combinations of these complete contractions.

The approximately invariant defining function $r=r^{\F}$ is too rough in
getting an expansion for $\psi$ analogous to that for $\varphi$, while
there is no hope of making $r$ exactly invariant. Instead, we consider a
family $\calF_M$ of defining functions of the germ $M$ of $\pa\Omega$ at a
point $p$ of reference such that $\calF_M$ is invariant under local
biholomorphic maps $\Phi\colon M\to\wt M$, that is, $r\in\calF_M$ if and
only if $\wt r\in\calF_{\wt M}$, where 
$\wt r\circ\Phi=|\det\Phi'|^{2/(n+1)}r$.  We also require that $\calF_M$ is
parametrized formally by $C^\infty(M)$. More precisely, $M$ is a formal
surface, $r$ is a formal function, and $C^\infty(M)$ should be replaced by
a space $C^\infty_{\rm formal}(M)$ of formal power series. If $M$ is
in normal form $N(A)$ with $p=0$, then\break 
$f\in C^\infty_{\rm formal}(M)$ is identified with the Taylor
coefficients $C=(C_{\alpha\c\beta}^l)$ of $f(z',\cz',v)$ as in
\eqref{mnf-rho}, so that the corresponding $r\in\calF_M$ has the
parametrization $r=r[A,C]$. Specific construction of $\calF_M$ is done by
lifting the \MA equation to $\C^*\times\Omega$ near $\C^*\times\pa\Omega$
and considering a family of local (or formal) asymptotic solutions, say
$\calF_M^{\rm aux}$, which is parametrized by $C^\infty_{{\rm 
formal}}(M)$. This is a refinement of Graham's construction \cite{G2} of
asymptotic solutions to the \MA equation in $\Omega$. Then, $\calF_M$
consists of the smooth parts of elements of $\calF_M^{\rm aux}$, and
the parametrization $C^\infty_{\rm formal}(M)\to\calF_M$ for $M=N(A)$
is given by the inverse map of $r\mapsto\pa^{n+2}_\rho r|_{\rho=0}$, which
comes from the parametrization of $\calF_M^{\rm aux}$.

Biholomorphic invariance of $\calF_M$ gives rise to an extension of the
$H$-action on the normal form coefficients $A$ to that on the pairs $(A,C)$.
In fact, a natural generalization of the CR invariant is obtained by
considering polynomials $I(A,C)$ in the variables $A_{\alpha\c\beta}^l$ and
$C_{\alpha\c\beta}^l$ such that
$$
  I(\wt A,\wt C)=|\det\Phi'(0)|^{-2w/(n+1)}I(A,C)
$$ 
as in \eqref{CRinv-transf}, for biholomorphic maps $\Phi$ and 
$(\wt A,\wt C)$ satisfying 
$r[\wt A,\wt C]\circ\Phi=|\det\Phi'|^{2/(n+1)}r[A,C]$. Such a polynomial
defines an assignment, to each pair $(M,r)$ with $r\in\calF_M$, of a
function $I[r]\in C^\infty(M)$:
\begin{equation}\label{general-CR-functional}
 I[r](p)=|\det\Phi'(p)|^{2w/(n+1)}I(A,C),
\end{equation} 
with $\Phi$ as in \eqref{CRinv-functional} and $(A,C)$ parametrizing 
$\wt r$ such that $\wt r\circ\Phi=|\det\Phi'|^{2/(n+1)}r$. We thus refer to
$I(A,C)$ as an {\it invariant of the pair $(M,r)$ of weight} $w$.

The problem for $\psi$ is then formulated as that of finding an asymptotic
expansion of $\psi$ in powers of $r\in\calF_{\pa\Omega}$ of the form
\begin{equation}\label{exp-psi-intro}
 \psi=\sum_{j=0}^\infty\psi_j[r]\,r^j+O^\infty(r)
 \quad \hbox{with} \ \ \psi_j[r]\in C^\infty(\overline\Omega),
\end{equation} 
such that each $\psi_j[r]|_{\pa\Omega}$ is an invariant of the pair
$(\pa\Omega,r)$ of weight $j+n+1$. As in the CR invariant case,  a class of
invariants of the pair $(\pa\Omega,r)$ is obtained by taking the boundary
value for linear combinations of complete contractions of tensor products
of covariant derivatives of the curvature of the metric $g[r]$. Elements of
this class are called {\it Weyl invariants}. We prove that all invariants
of the pair $(M,r)$ are Weyl invariants (see Theorems 4 and 5), so that the
expansion \eqref{exp-psi-intro} holds with $\psi_j[r]|_{\pa\Omega}$ given by
Weyl invariants of weight $j+n+1$ (see Theorem 1).

A CR invariant $I(A)$ is the same as an invariant of the pair $(M,r)$ which
is independent of the parameter $C$, so that $I(A)$ is a Weyl invariant
independent of $C$ (the converse also holds). That is, CR invariants are
the same as Weyl invariants independent of the parameter $C$ (see Theorem 2
which follows from Theorems 4 and 5). For Weyl invariants of low weight, it
is easy to examine the dependence on $C$. We have that all Weyl invariants
of weight $\leq n+2$ are independent of $C$ (see Theorem 3). This improves
the result of \cite{F3} and \cite{BEG} described above by weight $2$. If
$n=2$, we have a better estimate (see Theorem 3 again) which is consistent
with the results in \cite{HKN2}.

Introducing the parameter $C$ was inspired by the work of Graham \cite{G2}
on local determination of the asymptotic solution to the \MA equation in
$\Omega$. He proved approximate invariance, under local biholomorphic maps,
of the log term coefficients of the asymptotic solution, and gave a
construction of CR invariants of arbitrarily high weight. In our
terminology of Weyl invariants, these CR invariants are characterized as
complete contractions which contain the Ricci tensor of $g[r]$ (see Remark
\ref{rem-graham} for the precise statement).

This paper is organized as follows. In Section 1, we define the family
$\calF_M$ of defining functions and state our main results, Theorems 1, 2
and 3. Section 2 is devoted to the construction of the family $\calF_M$ and
the proof of its biholomorphic invariance. After reviewing the definition
of Moser's normal form, we reformulate, in Section 3, CR invariants and
invariants of the pair $(M,r)$ as polynomials in $(A,C)$ which are
invariant under the action of $H$. Then we relate these $H$-invariant
polynomials with those in the variables $R_{i\,\cj\, k\,\c l;ab\dots c}$ on
which $H$ acts tensorially, where $R_{i\,\cj\, k\,\c l;ab\dots c}$ are the
components of the curvature of $g[r]$ and its covariant derivatives. Using
this relation, we reduce our main Theorems 1--3 to the assertion that all
invariants of the pair $(M,r)$ are Weyl invariants. This assertion is
proved in two steps in Sections 4 and 5. In Section~4, we express all
invariants of the pair $(M,r)$ as $H$-invariant polynomials in 
$R_{i\,\cj\, k\,\c l;ab\dots c}$. In Section 5, we show that all such
$H$-invariant polynomials come from Weyl invariants, where invariant theory
of $H$ in \cite{BEG} is used essentially. In the final Section 6, we study
the dependence of Weyl invariants on the parameter $C$.

I am grateful to Professor Gen Komatsu, who introduced me to the analysis
of the Bergman kernel, for many discussions and encouragement along the way.

\section{Statement of the results}
 
\demo{{\rm 1.1.} Weyl functionals with exact transformation law}  
Our concern is a refinement of the ambient metric construction as in
\cite{F3}, \cite{BEG}.  Let $\Omega\subset\C^n$ be a smoothly bounded
strictly pseudoconvex domain and
$$
 J(u)=(-1)^n\det
  \left(
    \begin{array}{cc}
      u       & u_{\cj}\\
      u_i     & u_{i\,\cj}
     \end{array}
  \right)_{1\le i,\,j\le n}
\qtext{where}  u_{i\,\cj}=\pa_{z^i}\pa_{\c z^j}u.
$$ 
In \cite{F3}, \cite{BEG}, the construction started by choice of  a
defining function $r$, with $r>0$ in $\Omega$, satisfying
$J(r)=1+O^{n+1}(\pa\Omega)$, where $O^{n+1}(\pa\Omega)$ stands for a term
which is smoothly divisible by $r^{n+1}$. Such an  $r$ is unique  modulo
$O^{n+2}(\pa\Omega)$ and we denote the equivalence class by
$\calF^\F_{\pa\Omega}$. We here consider a subclass $\calF_{\pa\Omega}$ of
$\calF^\F_{\pa\Omega}$, which is  defined by lifting the (complex) \MA
equation (with Dirichlet boundary condition)
\begin{equation}\label{MA-J}
 J(u)=1 \hbox{ and } u>0\hbox{ in }\Omega,
 \quad u=0 \hbox{ on }\pa\Omega.
\end{equation} 
For a function $U(z^0,z)$ on $\C^*\times\c\Omega$, we set
$$
 J_\#(U)=(-1)^n\det\left(U_{i\,\cj}\right)_{0\le i,j\le n}
$$ 
and consider a \MA equation on $\C^*\times\Omega$:
\begin{equation}\label{MA}
 J_\#(U)=|z^0|^{2n}
  \hbox{ with } U>0 \hbox{ in }\C^*\times \Omega, \hbox{ and }
 \quad U=0 \hbox{ on }\C^*\times\pa\Omega. \qquad
\end{equation}  
If $U$ is written as $U(z^0,z)=|z^0|^2 u(z)$ with a function $u(z)$ on
$\Omega$, then \eqref{MA} is reduced to \eqref{MA-J} because
$J_\#(U)=|z^0|^{2n}J(u)$.  In \cite{G2}, Graham fixed
$r\in\calF_{\pa\Omega}^\F$ arbitrarily and constructed asymptotic solutions 
$u^\G$ to \eqref{MA-J} of the form
\begin{equation} \label{formal-solution-u}
  u^\G=r\sum_{k=0}^\infty\eta_k^\G\left(r^{n+1}\log r\right)^k
 \qtext{with} \eta_k^\G\in C^\infty(\c\Omega),
\end{equation} 
which are parametrized by the space $C^\infty(\pa\Omega)$ of initial data
(see Remark \ref{rem-asymptotic-solution} below).  Then
$U^\G=|z^0|^2u^\G$ are asymptotic solutions to \eqref{MA}.  We here modify
these asymptotic solutions and consider another class of asymptotic
solutions of the form
\begin{equation}\label{formal-solution}
  U=r_\#+r_\#\sum_{k=1}^\infty\eta_k\left(r^{n+1}\log r_\#\right)^k
 \qtext{with} \eta_k\in C^\infty(\c\Omega),
\end{equation}  
again parametrized by $C^\infty(\pa\Omega)$, where
$r_\#=|z^0|^2r$ with $r\in \calF_{\pa\Omega}^{\F}$. It should be emphasized
that $r$ is not prescribed but determined by $U$. Note also that $U$ is not
of the form $|z^0|^2u$ because $\log r_\#$  is not homogeneous in $z^0$. 
We call $r$ in \eqref{formal-solution} the {\it smooth part} of $U$ and
define $\calF_{\pa\Omega}$ to be the totality of the smooth parts of
asymptotic solutions to \eqref{MA} for $\pa\Omega$.

We identify two asymptotic solutions of the form \eqref{formal-solution} if
the corresponding functions $r$ and $\eta_k$ agree to infinite order along
$\pa\Omega$.  Then the unique existence of the asymptotic solution $U$ as
in \eqref{formal-solution} holds once the initial data are given in 
$C^\infty(\pa\Omega)$.

\specialnumber{1} \proclaim{Propostition} 
Let $X$ be a real vector field on $\c\Omega$ which is transversal to
$\pa\Omega${\rm .} Then for any $a\in C^\infty(\pa\Omega)${\rm ,} there exists a unique
asymptotic solution $U$ to {\rm \eqref{MA}} for $\pa\Omega$ such that the smooth
part $r$ satisfies  
\begin{equation}
  X^{n+2}r|_{\pa\Omega}=a.
  \label{initial-condition}
\end{equation}
\endproclaim

The lifted \MA equation \eqref{MA} and the asymptotic solutions of the form
\eqref{formal-solution} are introduced in order to obtain the following
exact transformation law for the smooth part $r$.

\specialnumber{2} \proclaim{Propostition}  
Let $\Phi\colon\Omega\to\wt\Omega$ be a biholomorphic map{\rm .}
Then  $r\in\calF_{\pa\Omega}$ if and only if 
$\wt r\in\calF_{\pa\wt\Omega}${\rm ,} where $\wt r$ is given by
\begin{equation}\label{trans-r}
 \wt r\circ\Phi=|\det\Phi'|^{2/(n+1)}r.
\end{equation} 
Here $\det\Phi'$ is the holomorphic Jacobian of $\Phi${\rm .}
\endproclaim

\numbereddemo{{R}emark}\label{rem-asymptotic-solution} 
For $u^\G$ in \eqref{formal-solution-u}, $\eta_0^\G=1+O^{n+1}(\pa\Omega)$
holds. To make $u^\G$ unique, Graham \cite{G2} used the boundary value of
$(\eta_0^\G-1)/r^{n+1}|_{\pa\Omega}$ as the initial data 
$a\in C^\infty(\pa\Omega)$, where $r$ is arbitrarily fixed.  It is also
possible to  make $u^\G$ unique by requiring $\eta_0^\G=1$ in 
\eqref{formal-solution-u}, in which case $r$ is determined by $u^\G$
(cf.~Lemma \ref{chose-r}).  Then we may write $r=r[u^\G]$ and consider the
totality of these, say $\calF_{\pa\Omega}^\G$.  However,
$\calF_{\pa\Omega}^\G$ does not satisfy the  transformation law
\eqref{trans-r} in Proposition 2; it is not the case that every $\wt
r=r[\wt u^\G]\in\calF_{\pa\wt\Omega}^\G$ is given by
\eqref{trans-r} with some $r=r[u^\G]\in\calF^\G_{\pa\Omega}$. Though the
proof requires some preparation  (cf.~Remark \ref{rem-proof}), this is
roughly seen by the fact that \eqref{trans-r} implies
$(\log \wt r)\circ\Phi=\log r+\log |\det\Phi'|^{2/(n+1)}$, which destroys
the condition $\wt\eta_0^\G=\eta_0^\G[\wt u^\G]=1$ (cf.~subsection
2.1).
\enddemo 

For each defining function $r\in\calF_{\pa\Omega}$, we define a
Lorentz-K\"ahler metric 
$$
  g[r]=\sum_{i,j=0}^n\frac{\pa^2\, r_\#}{\pa z^i\pa\c z^j}\,dz^i\,d\c z^j
  \qtext{on} \C^*\times\c\Omega \hbox{ near }\C^*\times\pa\Omega.
$$  
We call this metric $g=g[r]$ an {\em ambient metric} associated with
$\pa\Omega$.  From the ambient metric, we construct scalar functions as
follows. Let $R$ denote the curvature tensor of $g$ and  
$\R{p}{q}=\c\nabla{}^{q-2}\nabla^{p-2}R$ the successive covariant
derivatives, where $\nabla$ (resp.~$\c\nabla$) stands for the covariant
differentiation of type $(1,0)$ (resp.~$(0,1)$). Then a complete
contraction of the form
\begin{equation}\label{weyl-def}
  W_\#=\contr(\R{p_1}{q_1}\otimes\cdots\otimes\R{p_d}{q_d})
\end{equation}  
gives rise to a function $W_\#[r]$ on $\C^*\times\c\Omega$ near
$\C^*\times\pa\Omega$ once $r\in\calF_{\pa\Omega}$ is specified.  Here
contractions are taken with respect to the ambient metric for some pairing
of holomorphic and antiholomorphic indices. The {\em weight} of
$W_\#$ is defined by 
$
 w=-d+\sum_{j=1}^d (p_j+q_j)/2,
$ 
which is an integer because $\sum p_j=\sum q_j$ holds. 
By a {\it Weyl polynomial}, we mean a linear combination of $W_\#$ of the
form \eqref{weyl-def} of homogeneous weight. A Weyl polynomial gives a
functional for the pair $(\pa\Omega,r)$ which satisfies a  transformation
law under biholomorphic maps. To state this precisely, we make the
following definition.
 
\numbereddemo{{D}efinition}  
A Weyl polynomial $W_\#$  of weight $w$ assigns, to each pair
$(\pa\Omega,r)$ with $r\in\calF_{\pa\Omega}$, a function 
$W[r]=W_\#[r]|_{z^0=1}$ on $\c\Omega$ near $\pa\Omega$.  We call this
assignment $W\colon r\mapsto W[r]$   a {\it Weyl functional of weight} $w$
associated with~$W_\#$.
\enddemo

\specialnumber{3} \proclaim{Propostition}  
Let $W$ be a Weyl functional of weight $w${\rm .} Then{\rm ,} for $r$ and $\wt r$ as in
{\rm \eqref{trans-r},}
\begin{equation}\label{WF-inv-transf}
  W[\wt r]\circ\Phi=|\det\Phi'|^{-2w/(n+1)}W[r].
\end{equation}
\endproclaim  

We refer to the relation \eqref{WF-inv-transf} as a transformation law of
weight $w$ for $W$.

\numbereddemo{{R}emark}\rm 
Without change of the proof,  Proposition 1 can be localized near a boundary
point $p$. That is, we may replace $\pa\Omega$ by a germ $M$ of $\pa\Omega$
at $p$ or a formal surface, and $r$, $\eta_k$, $a$ by germs of smooth 
functions or formal power series about $p$. Then $\calF_{\pa\Omega}$ is a
sheaf $(\calF_{p,\c\Omega})_{p\in\pa\Omega}$.  Abusing notation, we write
$\calF_M$ in place of $\calF_{p,\c\Omega}$. Then Propositions 2 and 3 also
have localization, where $\Phi$ is a (formal) biholomorphic map such that
$\Phi(M)=\wt M$ with $\wt M$  associated to $\wt\Omega$. 
\enddemo

\demo{{\rm 1.2.} Invariant expansion of the Bergman kernel} For each
$r\in\calF_{\pa\Omega}$, we write the asymptotic expansion of the Bergman
kernel of $\Omega$ on the diagonal $K(z)=K(z,\c z)$ as follows:
\begin{equation}\label{bergman}
  K=\varphi[r]\,r^{-n-1}+\psi[r]\log r
  \qtext{with}
  \varphi[r],\psi[r]\in C^\infty(\c\Omega),
\end{equation} 
where we regard $\varphi=\varphi[r]$ and $\psi=\psi[r]$ as functionals  of
the pair $(\pa\Omega,r)$. Note that $\varphi[r]$ mod $O^{n+1}(\pa\Omega)$
and $\psi[r]$ mod $O^\infty(\pa\Omega)$ are independent of the choice of
$r$. In our first main theorem, we express these functionals in terms of
Weyl functionals.
\vglue6pt

\specialnumber{1} \proclaim{Theorem} 
For $n\ge2${\rm ,} there exist Weyl functionals $W_k$ of weight $k$ for
$k=0,1,2,\dots$ such that 
\begin{eqnarray}
\noalign{\vskip4pt}
  \varphi[r]& = &\ \sum_{k=0}^n  \ \, W_k[r]  \,r^k+O^{n+1}(\pa\Omega),
 \label{exp-phi}
\\ \noalign{\vskip4pt}
  \psi[r]& =&\sum_{k=0}^\infty W_{k+n+1}[r] \,r^k+O^\infty(\pa\Omega), \label{exp-psi}\\
\nonumber
\end{eqnarray} 
for any strictly pseudoconvex domain $\Omega\subset\C^n$ and any
$r\in\calF_{\pa\Omega}${\rm .} Here {\rm \eqref{exp-psi}} means that
$\psi[r]=\sum_{k=0}^m W_{k+n+1}[r]\,r^k+O^{m+1}(\pa\Omega)$ for any 
$m\ge 0${\rm .}
\endproclaim

The expansion \eqref{exp-phi} has been obtained in \cite{F3} and \cite{BEG},
where $r$ is any defining function satisfying $J(r)=1+O^{n+1}(\pa\Omega)$.
This condition is fulfilled by our $r\in\calF_{\pa\Omega}$.

\demo{{\rm 1.3. CR} invariants in terms of Weyl invariants} 
Suppose $\pa\Omega$ is in Moser's normal form (0.1) near $0$.  With the
real coordinates $(z',\c z',v,\rho)$,  we write the Taylor series about $0$
of $\pa^{n+2}_\rho r|_{\rho=0}$ for $r\in\calF_{\pa\Omega}$ as
\begin{equation}\label{pa-r}
 \pa^{n+2}_\rho r|_{\rho=0}=\sum_{|\alpha|, |\beta|, l\ge0}
  C_{\alpha\c\beta}^l z'{}^\alpha\c z'{}^\beta v^l.
\end{equation} 
Then for a Weyl functional $W$, the value $W[r](0)$ is expressed as a
universal polynomial $I_W(A,C)$ in the variables
$A_{\alpha\c\beta}^l, C_{\alpha\c\beta}^l$.  We call this polynomial a 
{\it Weyl invariant}\/ and say that $I_W$ is  {\it $\calC$-independent}\/ if
it is independent of the variables $C_{\alpha\c\beta}^l$. Our second main
theorem asserts that $\calC$-independent Weyl invariants give all CR
invariants.

\specialnumber{2} \proclaim{Theorem}\label{thm-weyl-general} 
All $\calC$\/{\rm -}\/independent Weyl invariants are {\rm CR} invariants{\rm ,} and vice versa{\rm .}
\endproclaim

It is not easy to determine which Weyl invariant $I_W$ is
$\calC$-independent when the weight $w$ of $I_W$ is high. If $w\le n+2$
(resp. $w\le5$) for $n\ge3$ (resp. $n=2$), then we can show that $W$ is
$\calC$-independent (Proposition \ref{prop-c-dependence}). Thus Theorem
\ref{thm-weyl-general} yields:

\specialnumber{3} \proclaim{Theorem} \label{thm-weyl} 
For weight $\le n+2${\rm ,} all Weyl invariants are {\rm CR} invariants and vice versa{\rm .}
Moreover{\rm ,} for $n=2${\rm ,} the same is true for weight $\le5${\rm .}
\endproclaim

In this theorem, the restriction on weight is optimal. In fact, there
exists a $\calC$-dependent Weyl invariant of weight $n+3$, or weight $6$
when $n=2$ (Proposition \ref{prop-c-dependence}). Thus, to obtain a
complete list of CR invariants for this or higher weights, one really needs
to select $\calC$-independent Weyl invariants.  This is a problem yet to be
studied.

\numbereddemo{{R}emark}\label{rem-weyl}\rm 
In the introduction, we  defined a Weyl invariant to be the boundary
value of a Weyl functional.  This definition is consistent with the one
given here as a polynomial $I_W(A,C)$.  In fact, $I_W(A,C)$ defines via
\eqref{general-CR-functional} an assignment, to each pair $(\pa\Omega,r)$,
of a function $I_W[r]\in C^\infty(\pa\Omega)$  which coincides with
$W[r]|_{\pa\Omega}$. This corresponds to the identification of a CR
invariant $I(A)$ with the boundary functional induced by $I(A)$.
\enddemo

\vglue-12truept
\section{Asymptotic solutions of the\\ complex  Monge-Amp\`ere equation}   

In this section we prove Propositions 1, 2 and 3. We first assume 
Proposition~1 and prove Propositions 2 and 3, the transformation laws of
$\calF_{\pa\Omega}$ and Weyl functionals.

\demo{{\rm 2.1.} Proof of Propositions {\rm 2} and {\rm 3}} 
For a biholomorphic map $\Phi\colon\Omega\to\wt\Omega$, we define the lift
$\Phi_\#\colon\C^*\times\Omega\to\C^*\times\wt\Omega$ by
\begin{equation}
  \Phi_\#(z^0,z)=\left(z^0\cdot[\det \Phi '(z)]^{{-1/(n+1)}},\Phi(z)\right),
  \label{Philift}
\end{equation}   
where a branch of $[\det \Phi ']^{{-1/(n+1)}}$ is arbitrarily chosen. Then
$\det\Phi'_\#(z^0,z)=[\det \Phi '(z)]^{{n/(n+1)}}$, so that
$$
 (|z^0|^{-2n}\det(\wt U_{i\,\cj}))\circ\Phi_\# =|z^0|^{-2n}\det((\wt
U\circ\Phi_\#)_{i\,\cj})
$$ 
for any function $\wt U$ on $\C^*\times\wt\Omega$. In particular, if $U$
is an asymptotic solution of \eqref{MA} for $\pa\Omega$, so is $\wt
U=U\circ\Phi^{-1}_\#$ for $\pa\wt\Omega$.  The expansion of $\wt U$ is
given by
$$
  \wt U=\wt r_\#+\wt r_\# 
  \sum_{k=1}^\infty \wt\eta_k\,({\wt r}^{\,n+1}\log\wt r_\#)^k,
$$  
where $\wt r\circ\Phi=|\det \Phi'|^{2/(n+1)}r$ and  
$\wt\eta_k\circ\Phi=|\det\Phi'|^{-2k}\eta_k$. It follows that $\wt r$ is the
smooth part of $\wt U$ if and only if 
$r=|\det\Phi'|^{-2/(n+1)}\wt r\circ\Phi$ is the smooth part of 
$U=\wt U\circ\Phi_\#$. This proves Proposition 2.
\medbreak

We next prove Proposition 3. Writing the transformation law \eqref{trans-r}
as $\wt r_\#\circ\Phi_\#=r_\#$ and applying $\pa\c\pa$ to it, we see that
$\Phi_\#\colon(\C^*\times\Omega,g[r])\to(\C^*\times\wt\Omega, g[\wt r])$ is
an isometry. If $W_\#$ is a Weyl polynomial of weight $w$, then
\begin{equation}
\label{invariance-W} W_\#[\wt r]\circ\Phi_\#=W_\#[r],
\end{equation}  
while the homogeneity of the ambient metric in $z^0$ implies
$$
 W_\#[r]=|z^0|^{-2w}W[r].
$$ 
Thus \eqref{invariance-W} is rewritten as \eqref{WF-inv-transf}, and 
Proposition 3 is proved.
\enddemo

\demo{{\rm 2.2.} Proof of Proposition {\rm 1}}  
We fix a defining function $\rho$ satisfying $J(\rho)=1+O^{n+1}(\pa\Omega)$
and introduce a nonlinear differential operator for functions $f$ on
$\C^*\times\Omega$:
$$
 \calM(f)=\det(U_{i\,\cj})/\det((\rho_\#)_{i\,\cj})
  \quad\hbox{with } U=\rho_\#\,(1+f).
$$  
Then $J_\#(U)=|z^0|^{2n}$  is written as 
\begin{equation}
  \calM(f)=J(\rho)^{-1}.
  \label{eqMf}
\end{equation}  
If $U$ is a series of the form \eqref{formal-solution}, then $f$ admits an
expansion 
$$
  f=\sum_{k=0}^\infty \eta_k (\rho^{n+1}\log\rho_\#)^k,
  \quad\hbox{where }
  \eta_k\in C^\infty(\c\Omega).
$$  
Denoting by $\cal A$ the space of all formal series of this form, we shall
construct solutions to \eqref{eqMf} in $\cal A$. 

We first study the degeneracy of the equation \eqref{eqMf} at the surface
$\C^*\times\pa\Omega$. Following \cite{G2}, we use a local frame
$Z_0,\dots,Z_n$ of $T^{(1,0)}(\C^*\times\c\Omega)$ near
$\C^*\times\pa\Omega$ satisfying:  
\begin{itemize}
\item[(1)] $Z_0=z^0(\pa/\pa z^0)$; 

\item[(2)] $Z_1,\dots,Z_{n-1}$ are orthonormal vector fields on
$\c\Omega$ with respect to the Levi form $\pa\c\pa\rho$ such that
$Z_j\rho=0$;

\item[(3)] $Z_n$ is a vector field on $\c\Omega$ such that 
 $Z_n\contraction\pa\c\pa\rho=\gamma\,\c\pa\rho$ for some 
$\gamma\in C^\infty(\c\Omega)$, $N\rho=1$ and $T\rho=0$, where $N=\Re Z_n$,
$T=\Im Z_n$.
\end{itemize}   
Using this frame, we introduce a ring $\calP_{\pa\Omega}$ of differential
operators on\break $\C^*\times\c\Omega$ that are written as polynomials of
$Z_0,\dots,Z_{n-1},\c Z_0,\dots,\c Z_{n-1}, T,\rho N$ with coefficients in
$C^\infty(\c\Omega,\C)$, the space of complex-valued smooth functions on
$\c\Omega$. In other words, $\calP_{\pa\Omega}$ is a ring generated by
$Z_0,\c Z_0$ and totally characteristic operators on $\c\Omega$ in the sense
of \cite{LM}. We first express $\calM$ as a nonlinear operator generated by
$\calP_{\pa\Omega}$. \enddemo

\proclaim{Lemma} \label{MAlin} 
Let $E=-(\rho N+1)(\rho N-2Z_0-n-1)${\rm .} Then{\rm ,}
\begin{equation}
  \calM(f)=1+Ef+\rho P_0f+Q(P_1f,\dots,P_l f)\qtext{for}f\in\calA,
  \label{Mfexp}
\end{equation}  
where $P_0,P_1,\dots,P_l\in\calP_{\pa\Omega}${\rm ,} and $Q$ is a polynomial
without constant and linear terms{\rm .}
\endproclaim 

\demo{Proof}    
Taking the dual frame $\omega^0,\dots,\omega^n$ of $Z_0,\dots,Z_n$, we set 
$\theta^j=z^0\omega^j$. Then, the conditions (1)--(3) imply $\theta^0=dz^0$,
$\theta^n=z^0\pa\rho$ and 
\begin{equation}
  \pa\c\pa\rho_\#=\rho\theta^0\wedge\c{\theta^0}
  +\theta^0\wedge\c{\theta^n}+\theta^n\wedge\c{\theta^0}
  -\sum_{i=1}^{n-1}\theta^i\wedge\c{\theta^i}+
  \gamma\theta^n\wedge\c{\theta^n}. \hskip.4in
  \label{pprho}
\end{equation}  
Using the coframe $\theta^0,\dots,\theta^n$, we define a Hermitian matrix
$A(f)=(A_{i\,\cj}(f))$ by  
$$
  \pa\c\pa\,\left(\rho_\#\,(1+f)\right)=\sum_{i,\, j=0}^n
  A_{i\,\cj}(f)\theta^i\wedge\c{\theta^j},
$$  
so that $\calM(f)=\det A(f)/\det A(0)$ holds. Let us compute $A(f)$. First,
$$
  \pa\c\pa\,\left(\rho_\#\,(1+f)\right)=(1+f)\pa\c\pa\,\rho_\#+
  \pa f\wedge\c\pa \rho_\#+ \pa \rho_\#\wedge\c\pa f
  +\rho_\#\pa\c\pa\,f.
$$  
For the first term on the right-hand side, we use \eqref{pprho}.  
The second and the third terms are respectively given by
$$
  \pa f\wedge\c\pa \rho_\#=\sum_{j=0}^{n}
  Z_j f\,\theta^j\wedge(\rho\,\c{\theta^0}+\c{\theta^n})
$$  
and its complex conjugate. Finally, for the last term,
$$
  \rho_\#\pa\c\pa f=\rho \sum_{i,j=0}^{n}
  (Z_i\c{Z_j}+E_{i\cj})f\,
   \theta^i\wedge\c{\theta^j}+\rho_\# Nf\pa\c\pa \rho,
$$  
where $E_{i\,\cj}\in\calP_{\pa\Omega}$ with $E_{0\cj}=E_{j\,\c 0}=0$ for
any $j$. Therefore, $A(f)$ modulo functions of the form $\rho Pf$,
$P\in{\cal P}_{\pa\Omega}$, is given by 
$$
  \left(
    \begin{array}{ccc}
      \rho&  0&  1+P_{0\c n}f\\
      0& -\delta_{i\cj}(1+f+\rho Nf)& *\\
      1+P_{n\c 0}f& *  & \gamma+P_{n\c n}f
    \end{array}
  \right),
$$ 
where $*$ stands for a function of the form $Pf, P\in\calP_{\pa\Omega}$, and
\begin{eqnarray*}
  P_{0\c n}&=&\c{P_{n\c 0}}=1+\rho \c{Z_n}+Z_0+\rho Z_0\c{Z_n},\\
  P_{n\c n}&=&\gamma+Z_n+\c{Z_n}+\rho Z_n\c{Z_n}+\gamma\rho N\\
           &=&\rho N^2+2N\quad\hbox{mod}\ \calP_{\pa\Omega}.
\end{eqnarray*}  Let $B(f)$ denote the matrix obtained from $A(f)$ by
dividing the first column by $\rho$ and multiplying the last row by $\rho$.
Then $B(f)$  modulo functions of the form $\rho Pf, P\in{\cal
P}_{\pa\Omega}$, is given by 
$$ 
  \left(\begin{array}{ccc}
    1 &  0       &  1+P_{0\c n}f\\
    {*} & -\delta_{i\cj}(1+f+\rho Nf) &*\\
    1+P_{n\c 0}f &0  & \gamma\rho+\rho^2N^2f+2\rho Nf
  \end{array}\right).
$$  
Noting that $\det A(f)=\det B(f)$, we get
\begin{eqnarray*}
  \calM(f)
  &=&
  1-\rho^2 N^2f-2\rho Nf+(n-1)(1+\rho N)f\\
  & &  +\ P_{n\c0}f+P_{0\c n}f+\rho P_0f+Q(P_1f,\dots,P_l f).
\end{eqnarray*}  
Using $Z_0f=\c{Z_0}f$, we obtain \eqref{Mfexp}.
\enddemo
\vglue6pt

To construct solutions to \eqref{eqMf} inductively, we introduce a
filtration 
$$
  \calA=\calA_0\supset\calA_1\supset\calA_2\supset\cdots,
$$  
where ${\cal A}_s$ denotes the space of all asymptotic series in $\calA$ of
the form

$$
  \rho^s\sum_{k=0}^\infty \alpha_k (\log\rho_\#)^k
  \qtext{with}\alpha_k\in C^\infty(\c\Omega).
$$  
\medbreak\noindent 
This filtration makes $\calA$ a filtered ring which is preserved by the
action of $\calP_{\pa\Omega}$. That is, $\calA_j\calA_k\subset\calA_{j+k}$
and $Pf\in\calA_j$ for each $(P,f)\in\calP_{\pa\Omega}\times\calA_j$. Hence
\eqref{Mfexp} yields $\calM(f+g)=\calM(f)+\calA_s$ for any $g\in\calA_s$. 
In particular, if $f\in\calA$ is a solution to the equation
$$
  \calM(f)=J(\rho)^{-1}\qtext{mod}\calA_{s+1},
  \leqno{\hbox{\eqref{eqMf}}{}_s}
$$  
\medbreak\noindent 
so is $\wt f=f+g$ for any $g\in\calA_{s+1}$. We shall show that this
equation admits a unique solution modulo $\calA_{s+1}$ if an initial
condition corresponding to \eqref{initial-condition} is imposed.

\vglue9pt

\proclaim{Lemma}\label{f-solution}  
{\rm (i)} An asymptotic series $f\in\calA$ satisfies {\rm \eqref{eqMf}}${}_n$ if
and only if $f\in\calA_{n+1}${\rm .}

{\rm (ii)}  Let $s\ge n+1${\rm .} Then{\rm ,} for any $a\in C^\infty(\pa\Omega)${\rm ,} the
equation {\rm \eqref{eqMf}}${}_s$ admits a solution $f_s${\rm ,} which is unique modulo
$\calA_{s+1}$ under the condition
\begin{equation}
  \eta_0=a\rho^{n+1}+O^{n+2}(\pa\Omega).
  \label{init-eta0}
\end{equation}
\endproclaim  

\vglue9pt

\demo{Proof}  
Since $f\in\calA$ satisfies $f=\eta_0$ mod $\calA_{n+1}$, it follows that
$\calM(f)=\calM(\eta_0)+\calA_{n+1}$.  Thus, recalling
$\calM(\eta_0)=J(\rho(1+\eta_0))/J(\rho)$, we see that \eqref{eqMf}${}_n$ is
reduced to 
$$
 J(\rho(1+\eta_0))=1+O^{n+1}(\pa\Omega).
$$  
\medbreak \noindent This is satisfied if and only if $\eta_0=O^{n+1}(\pa\Omega)$,  which is
equivalent to $f\in\calA_{n+1}$. Thus (i)  is proved.

To prove (ii), we first consider \eqref{eqMf}${}_{s}$ for $s=n+1$.  If
$f\in\calA_{n+1}$, then $\calM(f)=1+Ef+\calA_{n+2}$. Thus 
\eqref{eqMf}${}_{n+1}$ is equivalent to

\begin{equation}
  Ef=J(\rho)^{-1}-1 \qtext{mod}\calA_{n+2}.
  \label{MAn+1}
\end{equation}  
Writing $f=\rho^{n+1}(\alpha_0+\alpha_1\log\rho_\#)$ mod $\calA_{n+2}$, we
have $Ef=(n+2)\alpha_1\,\rho^{n+1}\break +\calA_{n+2}$. Hence, \eqref{MAn+1} holds
if and only if $(n+2)\alpha_1=(J(\rho)^{-1}-1)\rho^{-n-1}+O(\pa\Omega)$. 
Noting that $\alpha_0|_{\pa\Omega}$ is determined by \eqref{init-eta0}, we
get the unique existence of $f_{n+1}$ modulo $\calA_{n+2}$. 

For $s>n+1$, we construct $f_s$ by induction on $s$.  Assume that $f_{s-1}$
exists uniquely modulo $\calA_s$. Then we have
$\calM(f_{s-1}+g)=\calM(f_{s-1})+Eg+\calA_{s+1}$ for $g\in\calA_s$,  so
that \eqref{eqMf}${}_{s}$ is reduced to
\begin{equation}
  E[g]_s=[h]_s\qtext{with}h=J(\rho)^{-1}-\calM(f_{s-1})\in\calA_s,
  \label{linears}
\end{equation}  
where $[g]_s$ and $[h]_s$ stand for the cosets in $\calA_s/\calA_{s+1}$.
To solve this equation, we introduce a filtration of $\calA_s/\calA_{s+1}$:
$$
  \calA_s/\calA_{s+1}=\calA_s^{(l)}\supset
  \calA_s^{(l-1)}\supset\cdots\supset\calA_s^{(0)}
  \supset\calA_s^{(-1)}=\{0\},
$$   
where $l=[s/(n+1)]$ and 
$$
  \calA_s^{(t)}=\Bigl\{[g]_s\in\calA_s/\calA_{s+1}:
  g=\sum_{k=0}^t\eta_k(\rho^{n+1}\log\rho_\#)^k\in\calA_s\Bigr\}.
$$ 
Clearly, $\rho N\calA_s^{(t)}\!\!\subset\!\calA_s^{(t)}$ and
$Z_0\calA_s^{(t)}\!\!\subset\!\calA_s^{(t-1)}$.  Consequently, if we write
$[g]_s\!\in~\!\calA_s^{(m)}$ as 
$[g]_s=[\alpha_m\rho^s(\log\rho_\#)^m]_s+\calA_s^{(m-1)}$, then
$$
  E[g]_s=I(s)[\alpha_m\rho^s(\log\rho_\#)^m]_s+\calA_s^{(m-1)},
$$  
where $I(x)=-(x+1)(x-n-1)$. Thus, setting $F=1-I(s)^{-1}E$, we have 
$F[g]_s\in{\cal A}_s^{(m-1)}$ so that 
$F\calA_s^{(m)}\subset\calA_s^{(m-1)}$. In particular, $F^l=0$ on
$\calA_s^{(l)}$. Since $E=I(s)(1-F)$, the linear operator
$L=I(s)^{-1}\sum_{k=0}^{l-1} F^k$ satisfies $LE=EL=\hbox{id}$ on
$\calA_s^{(l)}$. Therefore, \eqref{linears} admits a unique solution
$[g]_s$, which gives a unique solution $f_s=f_{s-1}+g$ modulo $\calA_{s+1}$
of \eqref{eqMf}${}_s$.
\enddemo

The unique solution of \eqref{eqMf} with the condition \eqref{init-eta0} is
obtained by taking the limit of $f_s$ as $s\to\infty$.  More precisely, we
argue as follows.  For $a\in C^\infty(\pa\Omega)$,  we take a sequence
$\{f_s\}$ in Lemma \ref{f-solution}, and  write
$f_s=\sum\eta_k^{(s)}(\rho^{n+1}\log\rho_\#)^k$.  Then the uniqueness of
$f_s$ mod ${\cal A}_{s+1}$ yields
$\eta_k^{(s+1)}=\eta_k^{(s)}$ mod $O^{s-k(n+1)}(\pa\Omega)$.  This implies
the existence of $\eta_k\in C^\infty(\c\Omega)$ such that
$$\eta_k=\eta_k^{(s)} \hbox{ mod } O^{s-k(n+1)}(\pa\Omega)$$ for any $s$. Therefore,
the formal series $f=\sum_{k=0}^\infty\eta_k(\rho^{n+1}\log\rho_\#)^k$
satisfies $\calM(f)=J(\rho)^{-1}$ and \eqref{init-eta0}. The uniqueness
follows from that for each \eqref{eqMf}${}_{s}$.

We have constructed a solution $f\in \calA$ of \eqref{eqMf} and hence
obtained a formal series
\begin{equation}\label{Utemp}
  U=\rho_\#(1+f)=\rho_\#+\rho_\#\sum_{k=0}^\infty \eta_k
  (\rho^{n+1}\log\rho_\#)^k,
\end{equation}  
which solves \eqref{MA} to infinite order along $\C^*\times\pa\Omega$. 
In general, the series \eqref{Utemp} is not in the form
\eqref{formal-solution} because $\eta_0$ may not vanish. We next construct
a unique defining function $r$ such that $U$ is written in the form
\eqref{formal-solution}. In the following, we write $f=g$ mod
$O^\infty(\pa\Omega)$ if $f-g$ vanishes to infinite order along $\pa\Omega$.

\vglue6pt

\proclaim{Lemma}\label{chose-r} 
Let $f\in\calA_{n+1}${\rm .} Then there exists a unique defining function $r$ {\rm mod}
$O^\infty(\pa\Omega)$ such that $U=\rho_\#(1+f)$ is written in the form
{\rm \eqref{formal-solution}. }
\endproclaim  

\vglue6pt

\demo{Proof}  
Starting from $r_1=\rho$, we define a sequence of defining functions $r_s$,
$s=1,2,\dots$, by setting $r_{s+1}=r_s(1+\eta_{s,0})$,  where $\eta_{s,0}$
is the coefficient in the expansion
$U=r_{s\#}+r_{s\#}\sum_{k=0}^\infty\eta_{s,k}(r_s^{n+1}\log r_{s\#})^k$.
It then follows from $\log(r_{s+1\#})=\log(r_{s\#})+O^{s(n+1)}(\pa\Omega)$
that $\eta_{s,0}=O^{s(n+1)}(\pa\Omega)$, so that 
$r_{s+1}=r_s+O^{s(n+1)+1}(\pa\Omega)$.  We can then construct a defining
function $r$ such that $r=r_s$ mod $O^{s(n+1)+1}(\pa\Omega)$ for any $s$. 
With this $r$, the series $U$ is written as 
$U=r_{\#}+r_{\#}\sum_{k=1}^\infty\eta_k( r^{n+1}\log r_{\#})^k.$ 

Let us next prove the uniqueness of $r$. We take another defining function
$\wt r$ with the required property and write
$U=\wt r_{\#}\sum_{k=0}^\infty\wt\eta_k(\wt r^{n+1}\log\wt r_{\#})^k$. 
Setting $\phi= r/\wt r\in C^\infty(\c\Omega)$, we then have
$\wt\eta_0=\phi(1+\sum_{k=1}^\infty\eta_k(\rho^{n+1}\log\phi)^k)$. Since
$\wt\eta_0=1$,
\begin{equation}
  \frac{1}{\phi}=1+\sum_{k=1}^\infty\eta_k(r^{n+1}\log\phi)^k.
  \label{phi-relation}
\end{equation}  
This implies that if $\phi=1+O^m(\pa\Omega)$ then
$\phi=1+O^{m+n+1}(\pa\Omega)$. Therefore, $\phi=1$ mod
$O^\infty(\pa\Omega)$; that is, $\wt r=r$ mod $O^\infty(\pa\Omega)$.
\enddemo
 
\vglue6pt

We next examine the relation between the conditions
\eqref{initial-condition} and \eqref{init-eta0}. Writing $U=\rho_\#(1+f)$
in the form \eqref{formal-solution}, we have
$$
  r=\rho+\eta_0\,\rho+O^{2(n+1)}(\pa\Omega)
  =\rho+\rho^{n+2}a+O^{n+3}(\pa\Omega).
$$   
Applying $X^{n+2}$, we get 
$$
  X^{n+2}r=X^{n+2}\rho+(n+2)!\,(X\rho)^{n+2}a+O(\pa\Omega).
$$  
Since $X$ is transversal to $\pa\Omega$, that is, $X\rho|_{\pa\Omega}\ne 0$,
it follows that specifying $X^{n+2}r|_{\pa\Omega}$ is equivalent to
specifying $a$ in \eqref{init-eta0} when $\rho$ is prescribed. Therefore,
$f$ and thus $U=\rho_\#(1+f)$ are uniquely determined by the condition
\eqref{initial-condition}. This completes the proof of the first statement
of Proposition 1.

It remains to prove $r\in\calF^\F_{\pa\Omega}$; that is,
$J(r)=1+O^{n+1}(\pa\Omega)$. If we write $U=r_\#(1+f)$ then
$\calM(f)=J(r)^{-1}$, where $\calM$ is defined with respect to $\rho=r$. 
On the other hand, we have by Lemma \ref{f-solution}, (i) that
$f\in\calA_{n+1}$ and thus $\calM(f)=\calM(0)=1$ mod $\calA_{n+1}$.
Therefore, $J(r)^{-1}=1$ mod $\calA_{n+1}$; that is,
$J(r)=1+O^{n+1}(\pa\Omega)$.

\section{Reformulation of the main theorems}  

\demo{{\rm 3.1.} A group action characterizing {\rm CR} invariants}
We first recall the definition and basic properties of Moser's normal form
\cite{CM}. A real-analytic hypersurface $M\subset\C^n$ is said to be in
{\em normal form}\/ if it admits a defining function of the form
\begin{equation}
  \rho = 2u - |z'|^2-\sum
    _{|\alpha|, |\beta|\ge 2,\, l\ge 0}
  A_{\alpha\c\beta}^l\, {z'}^\alpha {\c z'}^\beta v^l,
  \label{mnf}
\end{equation}  
where the coefficients $(A_{\alpha\c\beta}^l)$ satisfy the following three
conditions:  (N1)  For each $p,\, q\ge 2$ and $l\ge0$, 
$\bA_{p\c q}^l=(A_{\alpha\c\beta}^l)_{|\alpha|=p,\, |\beta|=q}$ is a
bisymmetric tensor of type $(p, q)$ on $\C^{n-1}$;  (N2) 
$A_{\alpha\c\beta}^l=\c{A_{\beta\c\alpha}^l}$;  (N3)
$
  \tr\bA_{2\c2}^l=0,\ \tr^2\bA_{2\c3}^l=0,\ 
  \tr^3\bA_{3\c3}^l=0;
$ 
here $\tr$ denotes the usual tensorial trace with respect to
$\delta^{i\,\cj}$. We denote this surface by $N(A)$ with
$A=(A_{\alpha\c\beta}^l)$. In particular, $N(0)$ is the hyperquadric
$2\,u=|z'|^2$.

For any real-analytic strictly pseudoconvex surface $M$ and $p\in M$, there
exist holomorphic local coordinates near $p$  such that $M$ is in normal
form and $p=0$. Moreover, if $M$ is tangent to the hyperquadric to the
second order at $0$, a local coordinates change $S(z)=w$ for which $S(M)$ is
in normal form is unique under the normalization
\begin{equation} 
  S(0)=0,\quad S'(0)=\hbox{id},\quad \Im\frac{\pa^2 w^n}{\pa (z^n)^2}(0)=0.
  \label{normalization}
\end{equation} 
Even if $M$ is not real-analytic but merely $C^\infty$, there exists a
formal change of coordinates such that $M$ is given by a formal surface
$N(A)$, a surface which is defined by a formal power series of the form
\eqref{mnf}.  In this case, \eqref{normalization} uniquely determines $S$ 
as a formal power series.  We sometimes identify a formal surface $N(A)$ in
normal form with  the collection of coefficients $A=(A_{\alpha\c\beta}^l)$
and denote by $\cal N$ the real vector space of all $A$ satisfying (N1--3).

The conditions (N1--3) do not determine uniquely the normal form of a
surface: two different surfaces in normal form may be (formally)
biholomorphically equivalent. The equivalence classes of normal forms can
be written as orbits in $\calN$ of an action of the group of all fractional
linear transformations which preserve the hyperquadric and the origin. To
describe this action, let us first delineate the group explicitly.

In projective coordinates $(\zeta^0,\dots,\zeta^n)\in\C^{n+1}$ defined by
$z^j=\zeta^j/\zeta^0$, the hyperquadric is given by
$\zeta^0{\c\zeta^n}+\zeta^n{\c\zeta^0}-|\zeta'|^2=0$. Let $g_0$ denote the
matrix
\begin{equation}
  g_0=
    \left(
      \begin{array}{ccc}
        0&   0    &1\\
        0&-I_{n-1}&0\\
        1&   0    &0
      \end{array}
    \right).
  \label{g0matrix}
\end{equation}  
Then, the hyperquadric is preserved by   the induced linear fractional
transformation $\phi_h$ of $\C^n$ for $h\in\SU(g_0)$.  Clearly, $\phi_h$
leaves the origin $0\in\C^n$ fixed  if and only if $h$ is in the subgroup
$$
  H=\{ h\in \SU(g_0):\ h\,e_0=\lambda \,e_0,\ \lambda\in\C^*\},
$$  
where $e_0$ is the column vector ${}^t(1,0,\dots,0)$.

For $(h,A)\in H\times\calN$, the surface $\phi_h(N(A))$ always fits the
hyperquadric to the second order. Hence there is a unique transformation
$S$ normalized by \eqref{normalization} such that $S$ sends $\phi_h(N(A))$
to a surface $N(\wt A)= S\circ\phi_h(N(A))$ in normal form.  We set
$\Phi_{(h,A)}=S\circ\phi_h$ and $h.A=\wt A$.  Then the uniqueness of $S$
implies 
\begin{equation}
  \Phi_{(\wt h h,A)}=\Phi_{(\wt h,h.A)}\circ\Phi_{(h,A)}
  \qtext{for any}h,\wt h\in H.
  \label{transPhi}
\end{equation}  
Thus $H\times\calN\ni(h,A)\mapsto h.A\in\calN$ defines a left action of
$H$ on $\calN$. Moreover, any biholomorphic map $\Phi$ such that
$\Phi(0)=0$  and $\Phi(N(A))=N(\wt A)$ is written as $\Phi=\Phi_{(h,A)}$
for some $h$ satisfying $h.A=\wt A$. Therefore, the orbits of this
$H$-action are the equivalence classes of the normal form.

We now rewrite the definition of CR invariants in terms of this $H$-action.

\numbereddemo{{D}efinition}  
An $H$-{\it invariant of $\calN$ of weight $w$} is  a polynomial $I(A)$ in
the components of $A=(A_{\alpha\c\beta}^l)$ satisfying 
\begin{equation}
  h.I=\sigma_{w,w}(h)I\quad\hbox{for}\ h\in H,
  \label{invarianceP}
\end{equation} 
where $(h.I)(A)=I(h^{-1}.A)$ and $\sigma_{p,q}$ is the character of $H$
given by
$$
  \sigma_{p,q}(h)=\lambda^{-p}\c\lambda{}^{-q}
  \quad\hbox{when}\quad h\,e_0=\lambda\,e_0.
$$ 
We denote by $I^w(\calN)$ the vector space of all invariants  of $\calN$
of weight $w$. 
\enddemo

Let us observe that \eqref{invarianceP} is equivalent to
\eqref{CRinv-transf} and hence $I(A)$ is an\break $H$-invariant if and only if it
is a CR invariant. If $\Phi(N(A))=N(\wt A)$ and $\Phi(0)=0$, then there
exists an $h\in H$ such that $\wt A=h.A$ and $\Phi=\Phi_{(h,A)}$. Since
$\Phi'(0)=\phi_h'(0)$ and $\det\phi_h'(0)=\lambda^{-n-1}$,
\begin{equation} 
|\det\Phi'(0)|^{{-2w}/(n+1)}=|\det\phi'_h(0)|^{-2w/(n+1)}=\sigma_{-w,-w}(h).
  \label{det-sigma} 
\end{equation}  
Thus \eqref{CRinv-transf} is written as $I(h.A)=\sigma_{-w,-w}(h)I(A)$,
which is equivalent to \eqref{invarianceP}.

\demo{{\rm 3.2. } Action of $H$ on the defining functions}\label{H-action-to-r}
In the previous subsection, we  expressed the transformation law of CR
invariants in terms of the $H$-action on $\calN$. Proceeding similarly, let
us express the transformation law of defining functions \eqref{trans-r} by
using the group $H$.

We consider asymptotic solutions to \eqref{MA} for a surface $N(A)$ in
normal form, where $r$ and $\eta_k$ are regarded as real formal  power
series about $0\in\C^n$.  Let $\calF_{N(A)}$ denote the totality of the
smooth parts of asymptotic solutions for $N(A)$. Then, applying the
argument in the proof of Proposition 1, we can find for any real formal
power series 
\begin{equation}\label{h-c}
  h_C(z',\c z',v)=\sum_{|\alpha|,|\beta|,l\ge0}C_{\alpha\c\beta}^l
  z'{}^\alpha\c z'{}^\beta v^l,
\end{equation} 
a unique formal power series $r\in\calF_{N(A)}$ such that 
\begin{equation}\label{init}
  \pa^{n+2}_\rho r|_{\rho=0}=h_C,
\end{equation}  
where $\pa_\rho$ is the differentiation with respect to $\rho$ in the
(formal) coordinates $(z', \c z',v,\rho)$ given by \eqref{mnf}.  Thus,
denoting by $\calC$ the totality of the collections of coefficients
$C=(C_{\alpha\c\beta}^l)$ of \eqref{h-c}, we can define a map
$$
  \iota_1\colon\calN\times\calC\to\calF=\bigcup_{A\in\calN}\calF_{N(A)}
$$  
which assigns to each $(A,C)\in\calN\times\calC$ the element
$r\in\calF_{N(A)}$ satisfying \eqref{init}.

For $h\in H$ and $r\in\calF_{N(A)}$, we define $h.r=\wt r\in\calF_{N(h.A)}$
by
$$
 \wt r\circ\Phi_{(h,A)}=|\det\Phi_{(h,A)}'|^{2/(n+1)}r.
$$  
Then the map $H\times\calF\ni(h,r)\mapsto h.r\in\calF$ gives an
$H$-action on $\calF$ in virtue of \eqref{transPhi}, and it induces an
$H$-action on $\calN\times\calC$ via the bijection $\iota_1$. With respect
to this $H$-action, we can characterize  invariants of the pair $(M,r)$ as
$H$-invariant polynomials  of $\calN\times\calC$ defined as follows.

\numbereddemo{{D}efinition}  
An {\it $H$-invariant of $\calN\times\calC$ of weight $w$} is a polynomial
$I(A,C)$ in the variables $A_{\alpha\c\beta}^l,C_{\alpha\c\beta}^l$
satisfying \eqref{invarianceP} with $$(h.I)(A,C)=  I(h^{-1}.(A,C)).$$ 
The {\it totality} of such polynomials is denoted by $I^w(\calN\times\calC)$.
\enddemo

Observe that the projection $\calN\times\calC\to\calN$ is $H$-equivariant.
Hence an $H$-invariant of $\calN$ is regarded as an $H$-invariant of
$\calN\times\calC$, and $I^w(\calN)$ is identified with  a subspace of
$I^w(\calN\times\calC)$.

\demo{{\rm 3.3.} Tensorial realization} 
We next embed the $H$-space $\calF$, and also $\calN\times\calC$,  into a
tensor $H$-module by using the curvatures of the ambient metrics $g[r]$. 
Recall that, for $r\in\calF$, the ambient metric  is defined by the
K\"ahler potential $r_\#$, and $r_\#$  is a formal power series in $z^0,z$
(and $\c z^0,\c z$) about $e_0=(1,0)\in\C^*\times\C^n$. Since $r_\#$ is
homogeneous in $z^0$, it follows that  the ambient metric, the curvature
tensor $R_{j\,\c k\,l\,\c m}$ and the covariant derivatives  
$R_{j\,\c k\,l\,\c m,\gamma_1\cdots\gamma_p}$ are defined at each point
$\lambda\,e_0\in\C^*\times\C^n$, $\lambda\in\C^*$.

For simplicity of the notation, we write
$R^{(p,q)}=(R_{\alpha\c\beta})_{|\alpha|=p,|\beta|=q}$, where
$$
  R_{\alpha\c\beta}=
  \left\{
     \begin{array}{cl}
       R_{\alpha_1\c \beta_1\alpha_2\c \beta_2,
       \alpha_3\cdots\alpha_p\c \beta_3\cdots\c \beta_q}
         & \hbox{if}\ |\alpha|\ge2\hbox{ and }|\beta|\ge2;\\
       0 & \hbox{otherwise.}
     \end{array}
  \right.
$$  
Here, components of tensors are written with respect to the coordinates\break
$\zeta=(\zeta^0,\zeta^1,\dots,\zeta^n)\in\C^{n+1}$ given by
$$
\zeta^0=z^0,\ \zeta^1=z^0 z^1,\ \dots,\ \zeta^n=z^0 z^n.
$$ 
Then we have, at the point $e_0$,
\begin{equation}
\label{R-reduction0}
  \begin{array}{rl}
    R_{\alpha'0\alpha''\c\beta}
    &=(1-|\alpha'\alpha''|)R_{\alpha'\alpha''\c\beta},
    \\
    R_{\alpha\,\c{\beta'}\,\c0\,\c{\beta''}}
    &=(1-|\beta'\beta''|) R_{\alpha\,\c{\beta'}\,\c{\beta''}},
  \end{array}
\end{equation} 
for all lists $\alpha,\alpha',\alpha'',\beta,\beta',\beta''$ of indices in
$\{0,1,\dots,n\}$ with $|\alpha'|,|\beta'|>1$. This fact is a consequence
of the homogeneity of the K\"ahler potential $r_\#(\zeta,\c\zeta)$ in
$\zeta$ and $\c\zeta$; see the tensor restriction lemma in \cite{F3}.

We write down the transformation law of $R_{\alpha\c\beta}$ which comes
from  the\break $H$-action on $\calF$.  Let $W=\C^{n+1}$ denote the defining
representation of $\SU(g_0)$, hence also of $H$, by left multiplication on
the column vectors. We define\break $H$-modules 
$$
 \bT_s^{p,q}=(\otimes^{p,q} W^*)\otimes\sigma_{s-p,s-q},
\quad\hbox{for }p,q,s\in{\Bbb Z}\hbox{ with } p,q\ge0,
$$ 
where $\mathbold{\otimes}^{p,q} W^*=(\mathbold{\otimes}^p W^*)\otimes(\mathbold{\otimes}^q\c{W^*})$. We
denote by $\bT_s$ the $H$-submodule of $\prod_{p,q\ge 0}\bT_s^{p,q}$ 
consisting of all $T=(T_{\alpha\c\beta})_{|\alpha|,|\beta|\ge0}$ satisfying 
\addtocounter{equation}{1}
\def\reduction0{\theequation}
$$
  \begin{array}{rl}
    T_{\alpha'0\alpha''\c\beta}
    &=(s-|\alpha'\alpha''|)T_{\alpha'\alpha''\c\beta},
    \\
    T_{\alpha\,\c{\beta'}\,\c0\,\c{\beta''}}
    &=(s-|\beta'\beta''|) T_{\alpha\,\c{\beta'}\,\c{\beta''}},
  \end{array}
\eqno{\reduction0}_s
$$  
where $\alpha,\alpha',\alpha'',\beta,\beta',\beta''$ are lists of indices
with $|\alpha'|,|\beta'|> s$. Then, \eqref{R-reduction0} permits us to
define a map $\iota_2\colon{\calF}\to\bT_1$ by setting
$\iota_2(r)=(R_{\alpha\c\beta}|_{e_0})_{|\alpha|,|\beta|\ge0}$, where
$R_{\alpha\c\beta}$ are the components of the covariant  derivatives of the
curvature of~$g[r]$. \enddemo

\proclaim{Proposition} 
The map $\iota_2$ is $H$\/{\rm -}\/equivariant{\rm .} In particular{\rm ,} the image
$\calR=\iota_2(\calF)$ is an $H$\/{\rm -}\/invariant subset of\/ $\bT_1${\rm .  } 
\endproclaim

\demo{Proof}   
For $r\in\calF_{N(A)}$ and $h\in H$, set $\wt r=h.r$.  Then
$g[\wt r]=F_*g[r]$, where $F=(\Phi_{(h,A)})_\#$, so that
$$
 \wt R^{(p,q)}=F_*R^{(p,q)}\quad\hbox{for any $p,q\ge0$},
$$ 
where $R^{(p,q)}$ and $\wt R^{(p,q)}$ are curvatures of 
$g[r]$ and $g[\wt r]$, respectively.  Evaluating this formula at $e_0$, we have
$$
  \iota_2(h.r)=\left((F_*R^{(p,q)})|_{e_0}\right)_{p,q\ge0}.
$$  
Note that the right-hand side is independent of the choice of the lift
$F$.  We shall fix $F$ as in the next lemma and express 
$(F_*R^{(p,q)})|_{e_0}$ in terms of $R^{(p,q)}|_{e_0}$ and $h$. 

\proclaim{Lemma}\label{jacobim}  
There exists a lift $F$ of\/ $\Phi_{(h,A)}$ satisfying
$F(e_0)=\lambda\,e_0${\rm ,} where $\lambda=\sigma_{-1,0}(h)${\rm .} 
For such a lift $F${\rm ,} the Jacobian matrix $F'$ at $\lambda^{-1} e_0$ with
respect to $\zeta$ is $h${\rm .}
\endproclaim 

\demo{Proof of Lemma {\rm \ref{jacobim}}}  
We first note that the linear map $\zeta\mapsto h\,\zeta$ gives a lift of
$\phi_h$. This lift satisfies $(\phi_h)_\#(e_0)=\lambda\,e_0$ and
$(\phi_h)_\#'(\nu\,e_0)=h$ for any $\nu\in\C^*$. On the other hand, for a
map $S$ normalized by \eqref{normalization}, we can define $S_\#$ such that
$S_\#(e_0)=e_0$ and $S_\#'(e_0)=\id$; see Lemma N1 of \cite{F3}. Hence the
composition $F=S_\#\circ(\phi_h)_\#$ gives a lift of 
$\Phi_{(h,A)}=S\circ\phi_h$ satisfying $F(e_0)=\lambda\,e_0$. The Jacobian
matrix of $F$ at $\lambda^{-1}e_0$ is given by
$F'(\lambda^{-1} e_0)=S_\#'(e_0)\cdot(\phi_h)_\#'(\lambda^{-1}e_0)=h$.
\enddemo

By this lemma, we see that $(F_*R^{(p,q)})|_{e_0}$ is given by the usual
tensorial action of $h$ on 
$R^{(p,q)}|_{\lambda^{-1} e_0}\in\mathbold{\otimes}^{p,q} W^*$.  To compare
$R^{(p,q)}|_{\lambda^{-1} e_0}$ with $R^{(p,q)}|_{e_0}$,  we next consider
the scaling map $\phi(\zeta)=\lambda \zeta$.  Then, from the homogeneity of
$g$, we have $\phi_*g=|\lambda|^{-2}g$, so that
$\phi_*R^{(p,q)}=|\lambda|^{-2}R^{(p,q)}$. Thus we get
$R^{(p,q)}|_{\lambda^{-1}e_0}=
\lambda^{p-1}{\c\lambda}{}^{q-1}R^{(p,q)}|_{e_0}$. 
Summing up, we obtain
$(F_*R^{(p,q)})|_{e_0}=h.(R^{(p,q)}|_{e_0})\in\bT_1^{p,q}$. \phantom{it'scold here}
\hfill\qed
\medbreak

We have defined the following $H$-equivariant maps:
$$
  \calN\times\calC
  \stackrel{\iota_1}{\longrightarrow}\calF
  \stackrel{\iota_2}{\longrightarrow}\bT_1.
$$  
We set $\iota=\iota_2\circ\iota_1$. Inspecting the construction of
$\iota_1$ and $\iota_2$, we see that $\iota$ is a polynomial map in the
sense that each component of $\iota(A,C)=(T_{\alpha\c\beta}(A,C))$ is a
polynomial in the variables $(A_{\alpha\c\beta}^l,C_{\alpha\c\beta}^l)$.

We now define $H$-invariants of
$\calR=\iota(\calN\times\calC)\subset\bT_1$. 

\numbereddemo{{D}efinition}  
An $H$-{\it invariant of\, $\calR$ of weight $w$} is a polynomial $I(T)$ in
the components of $(T_{\alpha\c\beta})\in\bT_1$ satisfying 
$$
  I(h^{-1}.T)=\sigma_{w,w}(h)I(T)
  \quad\hbox{for any} \ (h,T)\in H\times\calR.
$$ 
Identifying two $H$-invariants which agree on $\calR$, we denote by
$I^w({\calR})$ the totality of the equivalence classes of all 
$H$-invariants of\/ $\calR$ of weight $w$.
\enddemo

This definition implies that $\iota$ induces an injection
\begin{equation}\label{iota-star}
  \iota^*\colon
  I^w(\calR)\ni I(T)\mapsto I(\iota(A,C))\in I^w(\calN\times\calC).
\end{equation} 
This map is also surjective by the following theorem.

\specialnumber{4} \proclaim{Theorem} 
There exists a polynomial map
$\tau\colon\bT_1\to\calN\times\calC$ such that $\tau\circ\iota=\id${\rm . } 
In particular{\rm ,} $\iota$ is injective and thus  the map {\rm \eqref{iota-star}} is
an isomorphism{\rm .}
\endproclaim

On the tensor space $\bT_1$, we can construct $H$-invariants by making
complete contractions of the form  
$$
  \contr\left(T^{(p_1,q_1)}\otimes\cdots\otimes T^{(p_d,q_d)}\right),
$$  
where $T^{(p,q)}\in\bT^{p,q}_1$ and the contraction is taken with respect to
the metric $g_0$. We call such $H$-invariants {\em elementary invariants}.
The next theorem asserts that $I^w(\calR)$ is spanned by the elementary
invariants of weight $w$.

\specialnumber{5} \proclaim{Theorem} 
Every $H$\/{\rm -}\/invariant of $\calR$ coincides on $\calR$ with a linear
combination of elementary invariants of\/ $\bT_1${\rm . }
\endproclaim

Once we know Theorems 4 and 5, we can easily prove the main theorems 
stated in Section 1.

\demo{{\rm 3.4.} Proofs of the main theorems using Theorems {\rm 4} and {\rm 5}} 
\enddemo

\demo{Proof of Theorem $1$} 
We here consider only the log term $\psi[r]$. The expansion of $\varphi[r]$
is obtained exactly in the same manner  if we note that $\varphi[r]$ is
defined only up to $O^{n+1}(\pa\Omega)$ and hence one should keep track of
the ambiguity in each step of the proof; see \cite{F3}.

We prove by induction on $m$ that there exist Weyl functionals $W_k$ of
weight $k$ such that
\def\exppsi{\theequation}
$$
  \psi[r]=\sum_{k=0}^{m-1} W_{k+n+1}[r]\, r^k+ J_m[r]\, r^m
  \qtext{for} J_{m}[r]\in C^\infty(\c\Omega).
  \leqno{\exppsi}_m
$$  
We interpret $(3.12)_0$ as $\psi[r]\in C^\infty(\c\Omega)$. Assuming
$(3.12)_m$, we seek a Weyl functional $W_{m+n+1}$ such that
$J_m[r]=W_{m+n+1}[r]$ on $\pa\Omega$ for any $\Omega$ and 
$r\in\calF_{\pa\Omega}$. Then, $(3.12)_{m+1}$ is obtained by setting 
$J_{m+1}[r]=(J_m[r]-W_{m+n+1}[r])/r$.

We first study the transformation law of $J_m$ under a biholomorphic map
$\Phi\colon\Omega\to\wt\Omega$. Let $\wt K$  be the Bergman kernel of
$\wt\Omega$ and $\psi[\wt r]$ its log term.  It then follows from $\wt
K\circ\Phi=|\det\Phi'|^{-2}K$ that
\begin{equation}\label{trans-psi}
  \psi[\wt r]\circ\Phi=|\det\Phi'|^{-2}\psi[r]\qtext{mod}
O^\infty(\pa\Omega).
\end{equation}  
Choosing $\wt r$ as in Proposition 2, \eqref{trans-r}, we obtain
\begin{equation}\label{trans-J}
  J_m[\wt r]\circ\Phi=|\det\Phi'|^{-2({m+n+1})/({n+1})}J_m[r]
  \quad\hbox{mod }O^\infty(\pa\Omega). \hskip.5in
\end{equation}  
That is, $J_m$ satisfies a transformation law of weight $m+n+1$. If
$\pa\Omega$ and $\pa\wt\Omega$ are locally written in normal form
$N(h^{-1}.A)$ and $N(A)$, respectively, then the restriction of
\eqref{trans-J} to $z=0$ gives
$$
  J_m[\wt r](0)=\sigma_{m+n+1,m+n+1}(h)J_m[r](0).
$$ 
{}From this formula, we can conclude 
$J_m[r](0)\in I^{m+n+1}(\calN\times\calC)$ if we know that $J_m[r](0)$ is a
polynomial in the components of $(A,C)\in\calN\times\calC$. Now we have:

\proclaim{Lemma}  
Assume that $\pa\Omega$ is locally written in normal form $N(A)${\rm ,}
$A\in\calN$ and that $r=\iota_1(A,C)${\rm .}   Then $J_m[r](0)$ is a universal
polynomial in the components of $(A,C)${\rm .}
\endproclaim 

\demo{Proof}  
Inspecting the proof of the existence of asymptotic expansion
\eqref{bergman} in \cite{F1} or \cite{BS}, we see that  the Taylor
coefficients of $\psi$ about $0$ are universal polynomials in $A$ (direct
methods of computing these universal polynomials are given in [B1,2]
and [HKN1,2]). On the other hand, by the constructions of
$r$ and Weyl invariants, we can show that the Taylor coefficients of $r$
and of $W_k[r]$ about $0$ are universal polynomials in $(A,C)$. Therefore,
we see by the relation $(3.12)_m$ that $J_m[r](0)$ is a universal
polynomial in $(A,C)$.
\enddemo

We now apply Theorems 3 and 4 to obtain a Weyl functional $W_{m+n+1}$ such
that $J_m[r](0)=W_{m+n+1}[r](0)$ for any $\pa\Omega$ in normal form and
any\break $r\in\calF_{\pa\Omega}$.  Noting that $J_m$ and $W_{m+n+1}$ satisfy the
same transformation law under biholomorphic maps, we conclude
$J_m[r]=W_{m+n+1}[r]$ on $\pa\Omega$ for any $\Omega$ by locally
transforming $\pa\Omega$ into normal form. We thus complete the  inductive
step.
\hfill\qed

\demo{Proof of Theorem $2$} 
Note that the isomorphism 
$\iota^*\colon I^w(\calR)\to I^w(\calN\times\calC)$ of Theorem 4 is given by
$W(R)\mapsto I_W(A,C)$. Thus Theorem 5 guarantees that each element of
$I^w(\calN\times\calC)$ is expressed as  a Weyl invariant $I_W(A,C)$ of
weight $w$. Noting that $I_W\in I^w(\calN)$ if and only if $I_W$ is
$\calC$-independent, we obtain Theorem 2. 
\enddemo 

\section{Proof of Theorem 4}
\demo{{\rm 4.1.} Reduction to finite-dimensional cases}  
We first write the map $\iota$ as a projective limit of maps $\iota_m$,
$m=1,2,\dots$, on finite-dimensional $H$-spaces.  Then, we can reduce the
proof of Theorem 4 to that of an analogous assertion for each
$\iota_m$ (Proposition 4.1 below).

To define the maps $\iota_m$, we introduce weights for the components of
$\calN$, $\calC$ and $\bT_s$ by considering the action of the matrix 
$$
  \delta_t=
  \left(
    \begin{array}{ccc}
      t&0      &0 \\
      0&I_{n-1}&0 \\
      0&0      &1/t
    \end{array}
  \right)
  \in H,\quad t>0.
$$  
The actions of $\delta_t$ on
$(A_{\alpha\c\beta}^l,C_{\alpha\c\beta}^l)\in\calN\times\calC$  and
$(T_{\alpha\c\beta})\in\bT_s$ are given by
\begin{eqnarray*}
  \delta_t.(A_{\alpha\c\beta}^l,C_{\alpha\c\beta}^l)
  &=&
  (t^{|\alpha|+|\beta|+2l-2}A_{\alpha\c\beta}^l,
  t^{|\alpha|+|\beta|+2l+2(n+1)}C_{\alpha\c\beta}^l),
  \\
  \delta_t.(T_{\alpha\c\beta})
  &=&
  (t^{\|\alpha\c\beta\|-2s}T_{\alpha\c\beta}).
\end{eqnarray*}  
Here $\|\alpha\c\beta\|$ is the strength of the indices defined by setting
$\| a_1\dots a_k\| =\| a_1\| +\cdots+\| a_k\|$ with 
$\|0\|=\|\c0\|=0$, $\| j\|=\|\cj\|=1$ for $1\le j\le n-1$  and $\| n\|=\|\c
n\|=2$.    The {\em weights} of the components $A_{\alpha\c\beta}^l$,
$C_{\alpha\c\beta}^l$ and $T_{\alpha\c\beta}$ are defined to be the halves
of degrees in $t$ with respect to the action of $\delta_t$:
$$
\frac{1}{2}(|\alpha|+|\beta|)+l-1,
\quad
\frac{1}{2}(|\alpha|+|\beta|)+l+n+1
\quad\hbox{and}\quad
\frac{1}{2}\|\alpha\c\beta\|-s,
$$
respectively.  We also define the weight for a monomial of
$(A_{\alpha\c\beta}^l,C_{\alpha\c\beta}^l)$ or $(T_{\alpha\c\beta})$ to be
the half of the degree with respect to $\delta_t$. Then the notion of
weight is consistent with that for $H$-invariants.

Let $[\calN]_m$ denote the vector space of all components
$A_{\alpha\c\beta}^l$ of weight $\le m$:
$$
  [\calN]_m=
  \left\{
    (A_{\alpha\c\beta}^l)_{(|\alpha|+|\beta|)/2+l-1\le m}:
    (A_{\alpha\c\beta}^l)\in\calN
  \right\}.
$$   
We also define $[\calC]_m$ and $[\bT_s]_m$ similarly. Then $[\calN]_m$,
$[\calC]_m$ and $[\bT_s]_m$ are finite-dimensional vector spaces such that
their projective limits as $m\to\infty$ are $\calN$, $\calC$ and $\bT_s$,
respectively.

Since $\iota$ is compatible with the action of $\delta_t$, each component
$T_{\alpha\c\beta}(A,C)$ of $\iota(A,C)$ is a homogeneous polynomial of
weight $||\alpha\c\beta||/2-1$. It follows that, if
$||\alpha\c\beta||/2-1\le m$, then $T_{\alpha\c\beta}(A,C)$ depends only on
the variables $(A,C)\in [{\cal N}]_m\times [{\cal C}]_m$, because all
components of $(A,C)\in{\cal N}\times{\cal C}$ have positive weight. 
We can thus define maps
$$
 \iota_m\colon[{\cal N}]_m\times[{\cal C}]_m\to[\bT_1]_m
\quad\hbox{by}\ \
  \iota_m(A,C)=(T_{\alpha\c\beta}(A,C))_{||\alpha\c\beta||/2-1\le m}.
$$ 
The projective limit of $\iota_m$ as $m\to\infty$ gives $\iota$. Therefore,
Theorem 4 is reduced to the following finite-dimensional proposition. \enddemo

\proclaim{Proposition}\label{prop-iota-m}  
\hskip-9pt There exist polynomial maps 
$\tau_m\colon[\bT_1]_m\to [\calN]_m  \times[\calC]_m${\rm ,} for $m=0,1,2\ldots\,${\rm ,}
such that $\tau_m\circ\iota_m=\id$ and
$\pi_m\circ\tau_{m}=\tau_{m-1}\circ\pi_{m}'$.  Here 
$\pi_m\colon[\calN]_m\times[\calC]_m\to[\calN]_{m-1}\times[\calC]_{m-1}$ 
and $\pi_{m}'\colon[\bT_1]_{m}\to[\bT_1]_{m-1}$  are the natural
surjections{\rm .}
\endproclaim

The projective limit of $\tau_m$ is a polynomial map 
$\tau\colon\bT_1\to\calN\times\calC$ such that $\tau\circ\iota=\id$. Thus
Theorem 4 follows from  Proposition
\ref{prop-iota-m}.

\demo{{\rm 4.2.} Linearization of the Monge-Amp\`ere equation} 
\label{Linearization-of-MA} 
We prove Proposition \ref{prop-iota-m} by using the inverse function
theorem. Our first task is to show the injectivity of
$d\iota_m|_0\colon T_0([\calN]_m\times[\calC]_m)\to[\bT_1]_m$,  the
differential of $\iota_m$ at $0=(0,0)$. Identifying
$T_0([\calN]_m\times[\calC]_m)$ with $[\calN]_m\oplus [\calC]_m$, we define
a linear map $d\iota|_0\colon\calN\oplus\calC\to\bT_1$ by the projective
limit of $d\iota_m|_0$.  Then we can prove the injectivity of each
$d\iota_m|_0$ by proving that of $d\iota|_0$.

Before starting the proof, we introduce some vector spaces of  formal power
series, which will be used in the rest of this paper. For $s\in{\Bbb Z}$,
let $\calE(s)$ denote the vector space of real formal power  series
$f(\zeta,\c\zeta)$ of $\zeta,\c\zeta$ about $e_0\in\C^{n+1}$  of homogeneous
degree $(s, s)$. Here we say that $f$ is  homogeneous of degree $(s, s')$
if $Zf=s\,f$ and $\c Zf=s'\,f$,  where $Z=\sum_{j=0}^n\zeta^j\pa_{\zeta^j}$.
The space $\calE(s)$ admits a natural $H$-action 
$(h.f)(\zeta,\c\zeta)=f(h^{-1}\zeta,\c{h^{-1}\zeta})$, and is embedded as
an $H$-submodule of $\bT_s$ by
$f\mapsto(T_{\alpha\c\beta})_{|\alpha|,|\beta|\ge0}$, where
$T_{\alpha\c\beta}=\pa_\zeta^\alpha\pa_{\c\zeta}^\beta f(e_0)$. Using this
expression, we define $H$-submodules of $\calE(s)$ for $s\ge0$:
\begin{eqnarray*}
  \calE_s &=&\{f\in\calE(s):
  T_{\alpha\c\beta}=0
  \hbox{ if }|\alpha|\le s\,\,\hbox{ or }\,\,|\beta|\le s\},
  \\
  \calE^s &=&\{f\in\calE(s):
  T_{\alpha\c\beta}=0
  \hbox{ if }|\alpha|> s\hbox{ and }|\beta|> s\}.
\end{eqnarray*} 
Then we have a decomposition of $\calE(s)$ as $H$-modules
\begin{equation}\label{calE-decomp}
  \calE(s)=\calE_s\oplus\calE^s.
\end{equation}

We next consider the restrictions of $f\in\calE(s)$ to the null cone 
$$
\calQ=\{\mu=\zeta^0\c\zeta^n+\zeta^n\c\zeta{}^0-|\zeta'|^2=0\}\subset\C^{n+
1}
$$ 
associated with $g_0$ and set
\begin{equation}\label{def-calJ}
  \calJ(s)=\{f|_\calQ:f\in\calE(s)\},\quad 
  \calJ^s =\{f|_\calQ:f\in\calE^s\}.
\end{equation} 
If we employ $(z^0,\c z^0,z',\c z',v)$ as coordinates of $\calQ$, then each
$f\in\calJ(s)$ is written as
\begin{equation}
  f=|z^0|^{2s}\sum_{|\alpha|,|\beta|,l\ge0}
  B_{\alpha\c\beta}^l z'{}^\alpha\c z'{}^\beta v^l.
 \label{relation-J}
\end{equation} 
Thus we may identify $\calJ(s)$ with the space of real formal power series
in $(z',\c z',v)$. Using this identification, we embed $\calN$ as a subspace
of $\calJ(1)$ and identify $\calC$ with $\calJ(-n-1)$,  so that the actions
of $\delta_t$ on $\calN$ and $\calC$ are compatible with those on $\calJ(1)$
and $\calJ(-n-1)$, respectively. It then follows from \cite{CM} that
\begin{equation}\label{moser-lem}
  \calJ(1)=\calN\oplus\calJ^1\quad 
  \hbox{(direct sum of vector spaces).}
\end{equation} 
This is the key equation in the proof of the uniqueness of the normal form.
\enddemo

\numbereddemo{{R}emark}  
(i) In the decomposition \eqref{moser-lem}, the vector space $\calJ^1$, as
well as $\calN$, is realized by a subspace of $\calJ(1)$.  In \cite{CM},
elements of $\calJ(1)$ are expressed by formal power series of
$(z',\c z',v)$, where $\calJ^1$ is identified via \eqref{relation-J} with
the range of a linear operator $L$ defined by
$$
  L(F)=\Re\left(\c z^1 f_1+\cdots+\c z^{n-1}f_{n-1}+f_n
          \right)\Big|_{u=|z'|^2/2}
$$ 
for $\C^n$-valued formal power series $F=(f_1,\dots,f_n)$ of $z$.

(ii) The $H$-action on $\calJ(1)$ is linear, whereas that on $\calN$
is nonlinear. These are defined differently, and $\calN$ in
\eqref{moser-lem} is not an $H$-invariant subspace of $\calJ(1)$.
Nevertheless, the tangent space $T_0\calN$ at the origin is
$H$-isomorphic to $\calJ(1)/\calJ^1$ as follows.  Let us tentatively
introduce an $H$-submodule $\calJ_3(1)$ of $\calJ(1)$ representing
surfaces close to those in normal form;  $\calJ_3(1)$ consists of
elements of $\calJ(1)$ vanishing to the third order at $e_0$, and each
$f\in\calJ_3(1)$ is identified with the surface $N(f)=N(B)$ via
\eqref{relation-J}. Then the $H$-action $(h,N(f))\mapsto\phi_h(N(f))$ on
$\calJ_3(1)$ is nonlinear, while the linearization coincides with the
action $\rho(h)$ on $\calJ(1)$. Now $\calN$ is a subset of $\calJ_3(1)$, and
the linearization of the $H$-action on $\calN$ is given by $(h,f)\mapsto
p(\rho(h)f)$, where $p\colon\calJ(1)\to\calN$ denotes the projection
associated with the decomposition \eqref{moser-lem}. This is the
$H$-action on $T_0\calN$, so that $T_0\calN$ is  $H$-isomorphic to the
quotient space $\calJ(1)/\calJ^1$.

(iii) In \cite{CM}, surfaces in normal form are constructed within
$\calJ_3(1)$ by induction on weight. It should be noted that the grading by
weight is different from the linearization.  More precisely, the lowest
weight part in the deviation from normal form is determined by the
projection $p\colon\calJ(1)\to\calN$, whereas higher weight parts are
affected nonlinearly by lower weight terms. We use a similar procedure in
the remaining part of this section, where the induction is implicit in
the inverse function theorem.
 \enddemo

Using these power series, we shall write down $d\iota |_0(A,C)$ explicitly. 
By definition,
$$
 d\iota |_0(A,C)=\frac{d}{d\e}\iota(\e\,A,\e\,C)\Big|_{\e=0}.
$$ 
To compute the right-hand side, we consider a family of asymptotic solutions
\begin{equation}\label{U-e}
  U_\e=r_{\e\#}+r_{\e\#}\sum_{k=1}^\infty\eta_{\e,k}
       \left(r_\e^{n+1}\log r_{\e\#}\right)^k
\end{equation} 
with $r_\e=\iota_1(\e\,A,\e\,C)$, and derive an equation such that $(d
U_\e/d\e)|_{\e=0}$  is a unique solution. 

\proclaim{Proposition}\label{proplinMA} 
{\rm (i)} 
$F=(d U_\e/d\e)|_{\e=0}$ admits an expansion of the form 
\begin{equation}
  F=\varphi+\eta\, \mu^{n+2}\log \mu,
  \quad\hbox{where }\varphi\in\calE(1),\ \eta\in\calE(-n-1).
  \label{expf}
\end{equation} 

{\rm (ii)} 
Let $\pa_{\mu}$ denote the differentiation with respect to $\mu$ for the
coordinates $(z^0,\c z^0,z',\c z',v,\mu)${\rm .} Then
\begin{eqnarray}
  \varphi|_{\calQ}=-f_A, 
  &\hbox{where}&
  f_A=|z^0|^2\sum A_{\alpha\c\beta}^l {z'}^\alpha{\c z'}^\beta v^l,
   \label{bcA}\\
  \pa_{\mu}^{n+2}\varphi|_{\calQ}=h_C, 
  &\hbox{where}&
  h_C=|z^0|^{-2n-2}\sum C_{\alpha\c\beta}^l {z'}^\alpha{\c z'}^\beta v^l.
  \label{bcC}
\end{eqnarray}

{\rm (iii)} Let 
$$
 \Delta=\pa_{\zeta^0}\pa_{\c\zeta^n}+\pa_{\zeta^n}\pa_{\c\zeta^0}
 -\sum_{j=1}^{n-1}\pa_{\zeta^j}\pa_{\c\zeta^j} 
$$ 
be the Laplacian for the ambient metric $g_0$ with potential $\mu${\rm .}  Then
$\Delta F=0.$ 

{\rm (iv)} 
For each $(f,h)\in\calJ(1)\oplus\calJ(-n-1)${\rm ,} there exists a unique function
$F$ of the form {\rm \eqref{expf}} satisfying $\Delta F=0$ and
\begin{equation}\label{initial-fg}
\varphi|_\calQ=-f,
 \quad
\pa_\mu^{n+2}\varphi|_\calQ=h.
\end{equation}
\endproclaim

In the proof of (i), we use the following lemma, which implies,
in particular, that the first variation of the higher log terms in
\eqref{U-e} vanishes at $\e=0$.

\proclaim{Lemma}
\label{e-dependence} Each log term coefficient of $U_\e$ satisfies
$\eta_{\e,k}=O(\e^k)${\rm .}
\endproclaim 

\demo{Proof}  
Setting $U_\e=r_{\e\#}(1+f_\e)$, we shall show that $f_\e\in O_\calA(\e)$,
where $O_\calA(\e)$ is the space of all formal series of the form
$$
  \sum_{k=1}^\infty\eta_{\e,k}\left(r_\e^{n+1}\log r_{\e\#}\right)^k
    \qtext{with} \eta_{\e,k}=O(\e^k).
$$  
Neglecting higher log terms in $f_\e$, we set 
$\wt f_\e=\eta_{\e,1}\,r_\e^{n+1}\log r_{\e\#}$. Then we have $\wt f_\e\in
O_\calA(\e)$ and $\calM_\e(\wt f_\e)=J(r_\e)^{-1}$ mod $\calA_{2(n+1)}^\e$,
where $\calM_\e$ and $\calA_{2(n+1)}^\e$ are defined  with respect to
$\rho=r_\e$. We can recover $f_\e$ from $\wt f_\e$ by the procedure of 
constructing asymptotic solutions used in the proof of Lemma
\ref{f-solution}.  This procedure consists of the operations of applying
$\calM_\e$ and solving the equation $E_\e[g_\e]_s=[h_\e]_s$ for $s>n+1$. 
These operations preserve the class $O_\calA(\e)$, so that
$\wt f_\e\in O_\calA(\e)$ implies $f_\e\in O_\calA(\e)$ as desired. 
\enddemo

\demo{Proof of Proposition $\ref{proplinMA}$} 
 (i) By Lemma \ref{e-dependence} above, we have
$$ 
  U_\e=r_{\e\#}+\mu^{n+2}|z^0|^{-2(n+1)}\eta_{\e,1}\log\mu+O(\e^2).
$$  
We thus get \eqref{expf} with 
$$
  \varphi=|z^0|^2(d\,r_\e/d\e)|_{\e=0}
  \quad\hbox{and}\quad
  \eta=|z^0|^{-2(n+1)}(d\,\eta_{\e,1}/d\e)|_{\e=0}.
$$ 

(ii) 
We now regard $f_A$ as an element of $\calE(1)$ which is independent of
$\mu$ in the coordinates $(z^0,\c z^0,z',\c z',v,\mu)$. Setting
$\mu_\e=\mu-\e\,f_A$, we then consider the Taylor expansion of $r_{\e\#}$
with respect to $\mu_\e$ in the coordinates 
$(z^0,\c z^0,z',\c z',v,\mu_\e)$:
\begin{equation}
  r_{\e\#}=\mu_\e+\sum_{k=1}^\infty
   a_{\e,k}(z^0,\c z^0,z',\c z',v)\,\mu_\e^k.
  \label{r-e-expand}
\end{equation} 
Since $r_{\e\#}=\mu+O(\e)$, we have $a_{\e,k}=O(\e)$.  Differentiating  
both sides of \eqref{r-e-expand} with respect to $\e$ and then setting
$\e=0$, we get
\begin{equation}
  \varphi=-f_A+\sum_{k=1}^\infty a'_k\,\mu^k,
  \qtext{where} a'_k=\left.\frac{d\,a_{\e,k}}{d\,\e}\right|_{\e=0}.
  \label{varphiexp}
\end{equation}  
Restricting this formula to $\mu=0$, we obtain \eqref{bcA}. We also have by
\eqref{varphiexp} that $\pa_\mu^{n+2}\varphi |_{\mu=0}=(n+2)!\,a'_{n+2}$. 
Therefore, noting that \eqref{init} is equivalent to
$a_{\e,n+2}=\e\,h_C/(n+2)!$, we obtain \eqref{bcC}.
\medskip

(iii) 
Recalling that $g_0$ is the flat metric given by the matrix
\eqref{g0matrix}, we have
$\Delta F=\tr(g_0^{-1}(\pa_{\zeta^i}\pa_{\c\zeta^j}F))$. Thus, by
$(\pa_{\zeta^i}\pa_{\c\zeta^j}U_\e)=
g_0+\e\,(\pa_{\zeta^i}\pa_{\c\zeta^j}F)+O(\e^2)$, we get
\begin{eqnarray*}
  \det( \pa_{\zeta^i}\pa_{\c\zeta^j}U_\e) 
  &=&\det\left(g_0\,(I_{n+1}+
       \e\, g_0^{-1}\,( \pa_{\zeta^i}\pa_{\c\zeta^j}F))\right)+O(\e^2)\\ 
  &=&\det g_0\,\det\left(I_{n+1}+
       \e\,g_0^{-1}\,(\pa_{\zeta^i}\pa_{\c\zeta^j}F)\right)+O(\e^2)\\
  &=&(-1)^n
\left(1+\e\,\tr(g_0^{-1}\,(\pa_{\zeta^i}\pa_{\c\zeta^j}F))\right)+O(\e^2)\\
  &=&(-1)^n (1+\e\,\Delta F) +O(\e^2).
\end{eqnarray*}  
Noting that $\det(\pa_{\zeta^i}\pa_{\c\zeta^j}U_\e)$ is independent of $\e$,
we obtain $\Delta F=0$. 
\medskip

(iv) By induction on $m$, we construct $F_m$ of the form \eqref{expf}
satisfying $\Delta F_m=O(\mu^m)$, where $O(\mu^m)$ stands for an expression
of the form\break $\mu^m(\varphi+\psi\log\mu)$ with $\varphi,\psi\in\calE(1-m)$. 
For $m=0$, we may take $F_0$ to be an arbitrary extension of $f$ to
$\calE(1)$.  Assume we have constructed $F_{m-1}$ with some $m>0$. When
$m<n+1$, we set $F_m=F_{m-1}+\varphi_m\,\mu^m$ with
$\varphi_m\in\calE(1-m)$, and get
\begin{eqnarray*}
  \Delta F_m
  &=&\Delta F_{m-1}+
      [\Delta,\mu^m]\varphi_m +\mu^{m}\Delta \varphi_m\\
  &=&\Delta F_{m-1} +m\,\mu^{m-1}(Z+\c Z+n+m)\varphi_m+O(\mu^m)\\
  &=&\Delta F_{m-1} +m(n+2-m)\mu ^{m-1}\varphi_m+O(\mu^m).
\end{eqnarray*} 
Thus $\Delta F_m=0$ holds with $\varphi_m=\mu^{1-m}\Delta F_{m-1}/m(n+2-m)$.
When $m\ge n+2$, we set $F_m=F_{m-1}+\mu^m(\varphi_m+\eta_m\log\mu)$ with
$\varphi_m,\eta_m\in\calE(1-m)$. Using 
\begin{equation}\label{commutator}
 \begin{array}{rl}
     [\Delta,\mu^m\log\mu]=
    &m\mu^{m-1}\log\mu\,(Z+\c Z+n+m)    \\
    &+\ \mu^{m-1}(Z+\c Z+2m+n),  
 \end{array}
\end{equation} 
we then obtain
\begin{eqnarray*}
  \Delta F_m
  &=&\Delta F_{m-1}
     +m(n+2-m)\mu^{m-1}\left(\varphi_m+\eta_m\log\mu\right)\\
  &  & +\ (n+2)\mu^{m-1}\eta_m+O(\mu^m).
\end{eqnarray*} 
If $m=n+2$, then 
$\Delta F_{n+2}=\Delta F_{n+1}+(n+2)\mu^{n+1}\eta_{n+2}+O(\mu^{n+2})$, 
so that $\varphi_{n+2}$ and $\eta_{n+2}$ are determined by 
\eqref{initial-fg} and $\Delta F_{n+2}=O(\mu^{n+2})$, respectively. 
If $m>n+2$, then $n+2-m\ne0$ and thus $\eta_m,\varphi_m$ are uniquely
determined by $\Delta F_m=O(\mu^m)$ as in  the case $m\le n+1$. This
completes the inductive step.  The solution to $\Delta F=0$ is then given by
the limit of $F_m$ as $m\to\infty$. The uniqueness assertion is clear by the
construction.
\enddemo

Using (iv), we define a linear map 
\begin{equation}\label{def-calF}
  \calL\colon\calJ(1)\oplus\calJ (-n-1)\to\calE(1),
\end{equation} 
where $\calL(f,h)=\varphi$ is the smooth part of the solution to
$\Delta(\varphi+\eta\mu^{n+2}\log\mu)=0$ satisfying \eqref{initial-fg}.
Setting $\varphi=\calL(f_A,h_C)$ for $(A,C)\in\calN\oplus\calC$, we get, by
(i) and (ii),
\begin{equation}
  r_{\e\#}={\mu}+\e\,\varphi+O(\e^2).
  \label{r-e-expansion}
\end{equation}  
Applying $\pa\c\pa$, we obtain
$$
  g_\e=g_0+\e\,\left(\varphi_{i\,\cj}\right)+O(\e^2).
$$  
Since $g_0$ is a flat metric, the curvature $R^\e_{\alpha\c\beta}$ of
$g_\e$ satisfies
\begin{equation}
  R^\e_{\alpha\c\beta}
   =\e\,\pa^\alpha_\zeta\pa^\beta_{\c\zeta}\,\varphi+O(\e^2)
  \qtext{for}|\alpha|,|\beta|\ge2,
  \label{R-linear}
\end{equation} 
where $\pa^{i\cdots j}_\zeta=\pa_{\zeta^i}\cdots\pa_{\zeta^j}$. Hence
$d\iota|_0(A,C)=(S_{\alpha\c\beta})$ is given by
$$
  S_{\alpha\c\beta}=
     \left\{ \begin{array}{ll}
       \pa^\alpha_\zeta\pa^\beta_{\c\zeta}\,\varphi(e_0)   & 
       \hbox{if $|\alpha|\ge2$ and $|\beta|\ge2$;}\\
       0 & \hbox{otherwise.}
     \end{array}\right.
$$  
Consequently, if we identify $(S_{\alpha\c\beta})$ with a series in
$\calE_1$, then 
\begin{equation}
  \label{composition}
 d\iota|_0(A,C)=\pi\circ \calL(f_A,h_C),
\end{equation} 
where $\pi\colon\calE(1)\to\calE_1$ is the projection with respect to the
decomposition \eqref{calE-decomp}. Using this expression, we now prove:

\proclaim{Proposition}\label{injectivity-of-diota} 
The map $d\iota|_0$ is injective{\rm . }
\endproclaim

{\it Proof}.
Assuming $d\iota|_0(A,C)=0$ for $(A,C)\in\calN\oplus\calC$, we shall prove
$(A,C)=(0,0)$. By \eqref{composition}, this assumption is equivalent to
$\calL(f_A,h_C)\in\calE^1$. We first write down the set
$(\Range\,\calL)\cap\calE^1$ explicitly.

\proclaim{Lemma}\label{lem-RangeL} 
{\rm (i)}
$\Range\,\calL=\wt\calH(1)${\rm ,} where
$$
 \wt\calH(1)=
  \left\{
    \varphi\in\calE(1):
    \Delta\varphi=c_n\,\mu^{n+1}\Delta^{n+2}\varphi,
   \ \Delta^{n+3}\varphi=0
  \right\}
$$ with $c_n=(-1)^{n+1}\left((n+1)!\right)^{-2}${\rm . }

{\rm (ii)} Let $\calH^1=\{\varphi\in\calE^1:\Delta\varphi=0\}${\rm .} Then
\begin{equation}
\label{range-L}
\wt\calH(1)\cap\calE^1=\calL(\calJ^1\oplus\{0\}) =\calH^1.
\end{equation}
\endproclaim 

\demo{Proof of Lemma $\ref{lem-RangeL}$} 
(i) For each $\varphi\in\Range\,\calL$, there exists $\eta\in\calE (-n-1)$ 
such that $F=\varphi+\eta\,\mu^{n+2}\log\mu$ satisfies $\Delta F=0$.
Using \eqref{commutator}, we then get
\begin{eqnarray*}
  \Delta F 
  &=&\Delta\varphi +\mu^{n+2}\log\mu\cdot \Delta\eta+
      [\Delta,\mu^{n+2}\log\mu]\,\eta\\
  &=&\Delta\varphi +(n+2)\mu ^{n+1}\eta+
      \mu ^{n+2}\log \mu\cdot\Delta\eta
\end{eqnarray*} 
so that $\Delta F=0$ is reduced to a system
\begin{equation}\label{system-varphi-eta}
  \left\{
    \begin{array}{c}
      \Delta\varphi+(n+2)\mu^{n+1}\eta=0,\\
      \Delta \eta=0.
    \end{array}
  \right.
\end{equation}  
Noting that this system implies
$\Delta^{n+2}\varphi=(-1)^n\,(n+1)!\,(n+2)!\eta$,  we replace $\eta$ in
\eqref{system-varphi-eta} by
$(-1)^{n}\Delta^{n+2}\varphi/[(n+1)!\,(n+2)!]$.  Then
\begin{eqnarray}
&&\label{eq-varphi}
  \left\{
    \begin{array}{c}
      \Delta\varphi=
        c_n\,\mu^{n+1}\Delta^{n+2}\varphi,\\
      \Delta^{n+3}\varphi=0.
    \end{array}
  \right.\\ 
\noalign{\hskip-.3in Thus we have $\Range\,\calL=\wt\calH(1)$.\hfill
\nonumber}
\end{eqnarray} 

\vglue-16pt
(ii) We first show that
$\wt\calH(1)\cap\calE^1\supset\calL(\calJ^1\oplus\{0\})$.
Since $\wt\calH(1)=\Range\,\calL$, it suffices to prove
$\calE^1\supset\calL(\calJ^1\oplus\{0\})$. For $f\in\calJ^1$, take its
extension $\wt f\in\calE^1$ such that $\pa_\mu\wt f=0$, and set
$$
 \varphi=\wt f-\mu\Delta\wt f/(n+1)\in\calE^1.
$$ 
Then, $\varphi|_\calQ=\wt f|_\calQ=f$, $\pa_\mu^{n+2}\varphi|_\calQ=0$ and
$$
 (n+1)\Delta\varphi=(n+1)\Delta\wt f-[\Delta,\mu]\Delta\wt
 f-\mu\Delta^2\wt f =-\mu\Delta^2\wt f=0,
$$
so that $\calL(f,0)=\varphi\in\calE^1$. We next show
$\wt\calH(1)\cap\calE^1\subset\calH^1$. For $\varphi\in\wt\calH(1)$, we have
$\Delta\varphi=c_n\mu^{n+2}\Delta^{n+2}\varphi$, while $\varphi\in\calE^1$
implies $\Delta^2\varphi=0$. Therefore, if
$\varphi\in\wt\calH(1)\cap\calE^1$, then $\Delta\varphi=0$ and thus
$\varphi\in\calH^1$. It remains to prove
$\calH^1\subset\calL(\calJ^1\oplus\{0\})$. But this is clear since each
$\varphi\in\calH^1$ satisfies 
$\varphi=\calL(\varphi|_\calQ,0)\break\in\calL(\calJ^1\oplus\{0\})$.
\enddemo

\vglue6pt 

We now return to the proof of Proposition \ref{injectivity-of-diota}. By
\eqref{range-L}, there exists an $f_1\in\calJ^1$ such that
$\calL(f_A-f_1,h_C)=0$.  The injectivity of $\calL$ then implies $f_A-f_1=0$
and $h_C=0$. Noting that \eqref{moser-lem} forces $f_A=f_1=0$, we get
$(A,C)=(0,0)$ as desired.
\hfill\qed

\vglue6pt 

\demo{{\rm 4.3.} Proof of Proposition {\rm \ref{prop-iota-m}}} 
By virtue of Proposition \ref{injectivity-of-diota}, we can apply the
inverse function theorem to $\iota_m$, and obtain a neighborhood $V$ of
$(0,0)\in [\calN]_m\times[\calC]_m$ such that $\iota_m|_V$ is an embedding. 
Noting that $\iota_m$ is compatible with the action of $\delta_t$, we see
that $\iota_m|_{V_t}$ is also an embedding, where
$V_t=\{\delta_t.(A,C):(A,C)\in V\}$. Since $t>0$ is arbitrary, it follows
that $\iota_m\colon [\calN]_m\times[\calC]_m\to [\bT_1]_m$ itself is an
embedding. Thus $[\calR]_m=\iota_m([\calN]_m\times[\calC]_m)$ is a
real-analytic submanifold of $[\bT_1]_m$, and there exists a real-analytic
inverse map $\iota_m^{-1}\colon[\calR]_m\to[\calN]_m\times[\calC]_m$.

We now construct $\tau_m$ inductively. The case $m=0$ is trivial because
$[\calN]_0\times[\calC]_0=\{(0,0)\}$. Assume we have gotten
$\tau_{m-1}(T)$. We construct $\tau_m(T)$ as follows. Denote the components
of $\iota_m^{-1}(R)$, $R\in [\calR]_m$, by 
$$
  (P_{\alpha\c\beta}^l(R),Q_{\alpha\c\beta}^l(R))
  \in[\calN]_m\times[\calC]_m.
$$ 
For $(P_{\alpha\c\beta}^l(R),Q_{\alpha\c\beta}^l(R))
\in[\calN]_{m-1}\times[\calC]_{m-1}$, define their polynomial extensions to
$[\bT_1]_m$ by the components of 
$\tau_{m-1}(T)=(P_{\alpha\c\beta}^l(T),Q_{\alpha\c\beta}^l(T))$. For the
components of weight $>m-1$,  we construct their polynomial extensions in
two steps. First, extend $P_{\alpha\c\beta}^l(R)$, $Q_{\alpha\c\beta}^l(R)$
to real-analytic functions $\wt P_{\alpha\c\beta}^l(T)$, 
$\wt Q_{\alpha\c\beta}^l(T)$ on $[\bT_1]_m$ in such a way that they have
homogeneous weight.  Next, neglect the monomials of degrees $>m$ in $\wt
P_{\alpha\c\beta}^l(T)$, $\wt Q_{\alpha\c\beta}^l(T)$ and define
polynomials $P_{\alpha\c\beta}^l(T)$, $Q_{\alpha\c\beta}^l(T)$.  These
polynomials are extensions of $P_{\alpha\c\beta}^l(R)$, 
$Q_{\alpha\c\beta}^l(R)$ because of the following lemma.

\proclaim{Lemma}\label{weight-of-R}  
A monomial $P(T)$ on $\bT_1$ vanishes on $\calR$ provided the weight is less
than the degree{\rm .}
\endproclaim  

\demo{Proof} Let $Q(A,C)=P(\iota(A,C))$ be of weight $w$.  Then, by the
assumption on $P(T)$, each monomial constituting $Q(A,C)$  has degree $>w$.
But such a monomial must be $0$ because all the variables
$A_{\alpha\c\beta}^l$ and $C_{\alpha\c\beta}^l$ have weight $\ge1$.  Thus
we have $Q(A,C)=0$, which is equivalent to $P(T)=0$ on $\calR$.
\enddemo
 
The collection $(P_{\alpha\c\beta}^l(T), Q_{\alpha\c\beta}^l(T))$ gives a
polynomial map $\tau_m(T)$ satisfying 
$\pi_m\circ\tau_m=\tau_{m-1}\circ\pi_m'$ and $\tau_m\circ\iota_m=\id$. This
completes the inductive step.

\numbereddemo{{R}emark}\rm\label{rem-proof}  
Using the method of linearization in this section, we can now prove the
statement in Remark \ref{rem-asymptotic-solution}. Given
$(A,C)\in\calN\oplus\calC$, let $u_\e^\G$ be the asymptotic solution  to
\eqref{MA-J} of the form \eqref{formal-solution-u} (with $\eta_1^\G=1$) for
the surface $N(\e A)$ satisfying \eqref{pa-r} with $\e C$ in place of $C$.
Then $F=|z^0|^2 (d u^\G_\e/d\e)|_{\e=0}$ can be written in the form 
\begin{equation}\label{G-lin}
  F=\varphi+\mu^{n+2}\eta\log(\mu|\zeta|^{-2}), 
  \qtext{where} \varphi\in\calE(1),\ \eta\in\calE(-n-1),\hskip.25in
\end{equation} 
which satisfies \eqref{bcA}, \eqref{bcC} and $\Delta F=0$. We denote $\varphi$
in \eqref{G-lin} by $\varphi[A,C]$, and set
$\calH^\G=\{\varphi[A,C]:(A,C)\in\calN\oplus\calC\}$.  Then, for
$\varphi,\wt\varphi\in\calH^\G$, $\varphi-\wt\varphi=O(\mu)$ if and only if
$\varphi-\wt\varphi=\mu^{n+2}\psi$ with $\psi\in\calE(-n-1)$ satisfying
$\Delta\psi=0$. Let $\calF_{\pa\Omega}^\G$ be the space of defining
functions in Remark \ref{rem-asymptotic-solution}, and assume that
$\calF_{\pa\Omega}^\G$  satisfies the transformation law
\eqref{trans-r}. Linearizing \eqref{trans-r}, we then see that, for each
$h\in\SU(g_0)$, 
\begin{equation}\label{transf-HG}
  \varphi\in \calH^\G\quad\hbox{if and only if}\quad
  \wt\varphi(\zeta):=\varphi(h\zeta)\in \calH^\G. 
\end{equation}
We next set $\wt F(\zeta)=F(h\zeta)$ for $F$ in \eqref{G-lin}. Then 
$\Delta\wt F=0$ and 
$\wt F=\widehat\varphi+\mu^{n+2}\wt\eta\log(\mu |\zeta^0|^{-2})$, where
\begin{equation}\label{whvarphi}
  \widehat\varphi(\zeta)=\wt\varphi(\zeta)+
  \mu^{n+2}\wt\eta\log(|\zeta^0/\wt\zeta^0|^2),\quad
  \wt\eta(\zeta)=\eta(\wt\zeta)\quad (\wt\zeta=h\zeta). \hskip.4in
\end{equation} 
Thus $\wt\varphi,\widehat\varphi\in\calH^\G$ whenever $\varphi\in\calH^G$.
It then follows from $\widehat\varphi-\wt\varphi=O(\mu)$ and
\eqref{whvarphi} that $\Delta\wt\eta\log(|\zeta^0/\wt\zeta^0|^2)=0$. But
this equation is not satisfied, e.g., by $h$ and $\varphi$ such that
$\wt\zeta^0=\zeta^0+i\zeta^n, \wt\zeta'=\zeta',\wt\zeta^n=\zeta^n$ and
$\varphi$ satisfying $\varphi=|\zeta^0|^2|\zeta^1/\zeta^0|^{2(n+2)}+O(\mu)$,
in which case $\eta=(-1)^{n+1}(n+2) |\zeta^0|^{-2(n+1)}$. This is a
contradiction, and we have proved the statement  in Remark
\ref{rem-asymptotic-solution}
\enddemo

\vglue-6truept
\section{Proof of Theorem 5} 
\vglue-6truept

{\rm 5.1.} {\it Linearization}.
We have seen in Theorem 4 that $\calR$ is a submanifold of $\bT_1$ with a
system of polynomial defining equations $\iota\circ\tau(T)-T=0$, where
$T=(T_{\alpha\c\beta})\in\bT_1$.   Using this fact, we first reduce the
study of $H$-invariants of $\calR$ to that of the invariants of the
$H$-module $T_0\calR$.  That is, we reduce Theorem 5 to the following:

\nonumproclaim{Theorem $5'$}  
Every $H$\/{\rm -}\/invariant of\, $T_0\calR$ is the restriction to\, $T_0\calR$ of a
linear combination of elementary invariants of\, $\bT_1$\/{\rm .}\/
\endproclaim

{\it Proof of Theorem $5$ using Theorem $5'$}.   We follow the argument of
\cite{F3}.   Taking an $H$-invariant $I$ of $\calR$ of weight $w$, we shall
show that, for any $N$, there exists a finite list of elementary invariants
$W_j$  such that $I$ is written in the form 
\begin{equation}
  I=\sum c_j W_j+Q_N\qtext{on} \calR \qtext{with}Q_N(T)=O(T^N),
  \label{Pinduction}
\end{equation}  
where $O(T^m)$ stands for a term (polynomial)  which does not contain
monomials of degree $<m$. Once this is proved, Theorem~5 follows. In fact, 
by taking $N$ so that $N>w$, we obtain by Lemma \ref{weight-of-R} that
$Q_N=0$ on $\calR$, that is, $I=\sum c_j W_j$ on $\calR$.

To prove \eqref{Pinduction}, we start by writing $I(T)=O(T^m)$ so that
$$
  I(T)=S^m(T)+O(T^{m+1}),
$$  
where $S^m$ is homogeneous of degree $m$.  Then $S^m$ is an $H$-invariant of
$T_0\calR $.  In fact, if we take a curve $\gamma_\e$ in $\calR$ such that
$\gamma_0=0$ and $(d\gamma_\e/d\e)|_{\e=0}=T\break\in T_0\calR $, then we have  
$S^m(T)=\lim_{\e\to 0}I(\gamma_\e)/\e^m$.  Since the right-hand side is
$H$-invariant, so is $S^m$ as claimed.  Therefore, we can find, by using
Theorem~$5'$, elementary invariants $W_j$  such that
\begin{equation}
  S^m=\sum c_j W_j+U,
  \label{Sm-tmp}
\end{equation} 
where $U$ is homogeneous of degree $m$ and vanishes on $T_0\calR $. We now
examine the remainder $U$. Let $\{P_i(T)\}_{i=1}^\infty$ be a system of
polynomials in the variables $T_{\alpha\c\beta}$ which defines $\calR$,
i.e., $\calR =\mathbold{\cap}_{i=1}^\infty\{P_i=0\}$, and let $p_i$ be the linear part
of $P_i$ so that $T_0\calR=\mathbold{\cap}_{i=1}^\infty\ker p_i$.  We write $U$ as a
finite sum $U=\sum U_i\,p_i$, where $U_i$ are homogeneous of degree $m-1$.
Then
$$
  U =\sum U_i\,(p_i-P_i)+\sum U_i P_i
    =\sum U_i\,(p_i-P_i)\qtext{on}\calR.
$$  
Noting $\sum U_i\,(p_i-P_i)=O(T^{m+1})$ and using \eqref{Sm-tmp}, we obtain
\eqref{Pinduction} for $N=m+1$. Repeating the same procedure for the
remainder $Q_{m+1}$, we obtain the expression \eqref{Pinduction}
inductively  for arbitrary $N$.
\hfill\qed

\demo{{\rm 5.2.} A short exact sequence characterizing $T_0\calR$} 
We further reduce Theorem $5'$ to an analogous theorem for a simpler
$H$-module of trace-free tensors. This is done by writing down  a system
of equations of $T_0\calR$ explicitly and giving a short  exact sequence
which characterizes $T_0\calR$, where 
$T_0\calR=d\iota|_0(\calN\oplus\calC)$ is regarded as a subspace of
$\calE_1$ as in subsection 4.2.

\proclaim{Proposition}\label{def-eq-T0R} 
{\rm (i)} The tangent space\/ $T_0\calR$ of\/ $\calR$ at\/ $0$ is given by
$$
 \wt\calH_1:= \wt\calH(1)\cap\calE_1=
  \left\{
    \varphi\in\calE_1:
    \Delta\varphi=c_n\,\mu^{n+1}\Delta^{n+2}\varphi,
   \ \Delta^{n+3}\varphi=0
  \right\}.
$$ 

{\rm (ii)} The following sequence is exact\/{\rm :}\/
\begin{equation}
   0\to\calH_1\hookrightarrow\wt\calH_1
   \stackrel{\Delta^{n+2}}{\longrightarrow}\calH(-n-1)\to 0,
   \label{short-exactsequence}
\end{equation} 
where $\calH(k)=\{\varphi\in\calE(k):\Delta\varphi=0\}$ and
$\calH_k=\calH(k)\cap\calE_k${\rm .}
\endproclaim 

{\it Proof}.  
(i) Since \eqref{range-L} implies $\pi\circ\calL(\calJ^1\oplus\{0\})=\{0\}$,
it follows from \eqref{moser-lem} that
$$
  \pi\circ\calL\left(\calN\oplus\calJ(-n-1)\right)=
  \pi\circ\calL\left(\calJ(1)\oplus\calJ(-n-1)\right).
$$ 
Using Lemma \ref{lem-RangeL} (i), we then get
\begin{eqnarray*}
  T_0\calR&=&\pi\circ\calL\left(\calN\oplus\calJ(-n-1)\right)\\
          &=&\pi(\Range\,\calL)\\ &=&\pi(\wt\calH(1))=\wt\calH_1.
\end{eqnarray*}

(ii) It is clear from the definition of $\wt\calH_1$ that 
$ 0\to\calH_1\hookrightarrow\wt\calH_1
 \stackrel{\Delta^{n+2}}{\longrightarrow}\calH(-n-1)$ is exact. It then
remains only to prove the surjectivity of $\Delta^{n+2}$. To show this, we
solve the equation 
\begin{equation}\label{Delta-eq}
  \Delta^{n+2}\varphi=\eta \quad\hbox{for $\eta\in\calH(-n-1)$ given}.
\end{equation} 
We need to find a solution $\varphi\in\wt\calH_1$.  But it suffices to find
$\varphi$ in $\wt\calH(1)=\wt\calH_1\oplus\calH^1$, because all
$\varphi\in\calH^1$ satisfy $\Delta^{n+2}\varphi=0$.  Next, we
follow the argument of \cite[Prop.~4.5]{EG2}.  We first recall a lemma
in \cite{EG1}.

\proclaim{Lemma}
 For $k\le0${\rm ,} the restriction $\calH(k)\ni\eta\mapsto
\eta|_\calQ\in\calJ(k)$ is an isomorphism{\rm .} 
\endproclaim  

\demo{Proof} 
This amounts to proving the unique existence of a solution $\eta\in\calH(k)$
to the equation $\Delta\eta=0$ under the condition
$\eta|_\calQ=f\in\calJ(k)$.  The proof is a straightforward modification of
that of our Proposition \ref{proplinMA}, (iv).
\enddemo

By this lemma, we can reduce \eqref{Delta-eq} to an equation for
$f\in\calJ(1)$:
\begin{equation}\label{Delta-eq-2}
  \Delta^{n+2}\calL(f,0)|_\calQ=g,
  \quad\hbox{where $g=\eta|_\calQ\in\calJ(-n-1)$}.
\end{equation} 
We write down the left-hand side  with the real coordinates
$(t,x)=(t,x^1,\dots,$ $x^{2n})$ of
$\calQ$, where $2\zeta^0=t+i\,x^1$, $2\zeta^j=x^{2j}+i\,x^{2j+1}$,
$j=1,\dots,n-1$, and $x^{2n}=2\,\Im\zeta^n$. 

\proclaim{Lemma} 
Let  $\varphi$ be a formal power series about $e_0\in\C^{n+1}$ of
homogeneous degree $2$ in the sense that $(Z+\c Z)\varphi=2\varphi${\rm .} Then 
\begin{equation}\label{Delta-Delta}
  (\Delta^{n+2}\varphi)|_\calQ=\Delta_x^{n+2}(\varphi |_\calQ),
  \quad\hbox{where }
  \Delta_x=2\pa_{x^1}\pa_{x^{2n}}-\sum_{j=2}^{2n-1}\pa_{x^j}^2. \hskip.5in
\end{equation}
\endproclaim 

\demo{Proof}  In terms of the coordinates $(t,x,\mu)$, the Laplacian
$\Delta$ is written as
$$ 
  \Delta=\Delta_x+
  \left(
     \mu\pa_\mu+E +n+1
  \right)\pa_\mu,
  \quad\hbox{where }E=t\pa_t+\sum_{j=1}^{2n}x^j\pa_{x^j}.
$$ 
Writing $\varphi(t,x,\mu)=\varphi_0(t,x)+\mu\,\psi(t,x,\mu)$, we have
$$\Delta^{n+2}\varphi=\Delta^{n+2}_x\varphi_0+\Delta^{n+2}(\mu\,\psi).$$ 
Noting that $\psi$ is homogeneous of degree $0$, we have
\begin{eqnarray*}
  \Delta^{n+2}(\mu\,\psi)
  &=&[\Delta^{n+2},\mu]\psi+O(\mu)\\ \noalign{\vskip4pt}
  &=&(n+2)(Z+\c Z+2n+2)\Delta^{n+1}\psi+O(\mu)\\ \noalign{\vskip4pt}
  &=&O(\mu).
\end{eqnarray*} 
Therefore, $\Delta^{n+2}\varphi=\Delta^{n+2}_x\varphi_0+O(\mu)$, which is
equivalent to \eqref{Delta-Delta}.
\enddemo

\phantom{space}

Since \eqref{Delta-Delta} implies
$\Delta^{n+2}\calL(f,0)|_\calQ=\Delta_x^{n+2}f$, we can reduce
\eqref{Delta-eq-2} to
\begin{equation}\label{equation-x}
  \Delta_x^{n+2}f=g.
\end{equation} 
It is a standard fact of harmonic polynomials that, for each polynomial
$q(x)$ of homogeneous degree $k$, there exists a polynomial $p(x)$ of
homogeneous degree $k$ such that
$$
  \Delta_x^{n+2}\mu_x^{n+2}p=q,\quad\hbox{where } 
  \mu_x=2 x^1 x^{2n}-\sum_{j=2}^{2n-1}(x^j)^2.
$$
We apply this fact to solving \eqref{equation-x}. Writing 
$g(t,x)=\sum_{j=2n+2}^\infty q_j(x)\,t^{-j}$ with polynomials $q_j(x)$ of
homogeneous degree $j-2n-2$, we take, for each $j$, a polynomial $p_j(x)$
such that $\Delta_x^{n+2}\mu_x^{n+2}p_j=q_j$. Then 
$$
  \Delta^{n+2}_x\wt f=g,\quad\hbox{where }
  \wt f=\mu_x^{n+2}\sum_{j=2n+2}^\infty p_j(x)\,t^{-j}.
$$ 
It is clear that $\wt f$ is homogeneous of degree $2$, though $\wt f$
may not be contained in $\calJ(1)$.  Let us write
$\wt f=\sum_{p+q=2}f^{(p,q)}$, where $f^{(p,q)}$ is homogeneous of degree 
$(p,q)$.  We then see by \eqref{Delta-Delta} that $\Delta_x^{n+2}f^{(p,q)}$
is homogeneous of degree\break $(p-n-2,q-n-2)$. Setting $f=f^{(1,1)}\in\calJ(1)$,
we thus obtain \eqref{equation-x}. The proof of Proposition \ref{def-eq-T0R}
is complete.
\hfill\qed\medbreak

 \phantom{space}

Now we use the exact sequence \eqref{short-exactsequence} and reduce
Theorem $5'$ to: \eject
\nonumproclaim{Theorem $5''$} 
Every $H$-invariant of $\calH_1\oplus\calH (-n-1)$ is realized by  the
restriction to $\calH_1\oplus\calH (-n-1)$ of a linear combination of
elementary invariants of $\calE_1\oplus\calE(-n-1)${\rm .} Here{\rm ,} an elementary
invariant of $\calE_1\oplus\calE(-n-1)$ is defined to be a complete
contraction of the form  
$$
  \contr
  \left(
    R^{(p_1,q_1)}\otimes\cdots\otimes R^{(p_d,q_d)}\otimes
    E^{(p'_1,q'_1)}\otimes\cdots\otimes E^{(p'_{d'},q'_{d'})} 
  \right),
$$  
where $R^{(p,q)}=(R_{\alpha\c\beta})_{|\alpha|=p,|\beta|=q}$ and 
$E^{(p,q)}=(E_{\alpha\c\beta})_{|\alpha|=p,|\beta|=q}$ with
$(R_{\alpha\c\beta},E_{\alpha\c\beta})\break\in\calE_1\oplus\calE(-n-1)
\subset\bT_1\oplus\bT_{-n-1}${\rm .}
\endproclaim

{\it Proof of Theorem $5'$ using Theorem $5''$}. 
 We embed $T_0\calR$ into $\calE_1\oplus\calE(-n-1)$ as a subspace 
$\calH =\{(R,\Delta^{n+2}R):R\in T_0\calR\}$ by identifying
$R=(R_{\alpha\c\beta})$ with a formal power series in $\calE_1$.  We wish
to find, for any $N$,  a list of elementary invariants $\{W_j\}$ of $\bT_1$
such that $I$ is written in the form
\begin{equation}
  I=\sum c_j W_j+O(E^N)\qtext{on} \calH,
  \label{Pexp-tmp}
\end{equation} 
where $O(E^N)$ is a polynomial in $(R_{\alpha\c\beta},E_{\alpha\c\beta})$
consisting of monomials of degree $\ge N$ in $E$. If $(n+1)N$ is greater
than the weight of $I$, then the error term $O(E^N)$ vanishes, because each
component $E_{\alpha\c\beta}$ has weight $\ge n+1$, and the reduction to
showing \eqref{Pexp-tmp} is done. 

The proof of \eqref{Pexp-tmp} goes as an analogy of the procedure of
linearization. Writing 
$$
 	I(R,E)=S^m(R,E)+O(E^{m+1}),
$$   
where $S^m$ is homogeneous of degree $m$ in $E$, we show that $S^m$ is an
$H$-invariant of $\calH_1\oplus\calH (-n-1)$.  For any
$(R,E)\in\calH_1\oplus\calH (-n-1)$, we use \eqref{short-exactsequence} and
choose $\wt R$ such that $(\wt R,E)\in \calH $. Then we have
$S^m(R,E)=\lim_{\e\to0}I(R+\e\wt R,\e E)/\e^m$. Since the right-hand side is
$H$-invariant, so is $S^m$ as claimed.  Therefore we can find, by Theorem
$5''$, a list of elementary invariants $W_j(R,E)$ of 
$\calE_1\oplus\calE(-n-1)$ and write $S^m$ as
$$
  S^m=\sum c_j W_j+U,
$$  
where $U$ vanishes on $\calH_1\oplus\calH (-n-1)$ and is homogeneous of
degree $m$ in $E$. Note that each elementary invariant $W_j(R,E)$ coincides
on $\calH$ with the elementary invariant $W_j(R,\Delta^{n+2}R)$ of
$\bT_1$. We next study the remainder $U$. Recall by
Proposition~\ref{def-eq-T0R} that $\calH $ is written as
$$
  \calH =\{(R,E): 
  P_i(R)=\wt Q_i(E), Q_i(E)=0\hbox{ and } E=\Delta^{n+2}R\},
$$
 where $\{P_i(R)\}_{i=1}^\infty$ and $\{Q_i(E), \wt Q_i(E)\}_{i=1}^\infty$
are systems of linear functions on $\calE_1$ and $\calE(-n-1)$,
respectively, such that 
$$
  \calH_1=\cap_i\ker P_i
  \quad\hbox{and}\quad
  \calH(-n-1)=\cap_i\ker Q_i.
$$ 
Using the defining functions $\{P_i,Q_i\}$ of $\calH_1\oplus\calH(-n-1)$, 
we can express $U$ as a finite sum $U=\sum U_i P_i+\sum V_i Q_i$, where
$U_i$ (resp.~$V_i$) are homogeneous of degree $m$ (resp.~$m-1$) in $E$. We
thus write $U$ in the form
\begin{equation}
  U=\sum U_i\,\wt Q_i+\sum U_i(P_i-\wt Q_i)+\sum V_i Q_i
  \label{Uexpression}
\end{equation}  
and find that $U=\sum U_i\wt Q_i=O(E^{m+1})$ on $\calH$, because the last
two sums in \eqref{Uexpression} vanish on $\calH$. We have shown
\eqref{Pexp-tmp} for $N=m+1$. Repeating the same procedure for the
remainder, we obtain \eqref{Pexp-tmp} for arbitrary $N$.
\hfill\qed

\demo{{\rm 5.3.} Proof of Theorem $5''$}  
Since $H$ acts on $\calH_1\oplus\calH(-n-1)$ by linear transformations, we
may restrict our attention to the $H$-invariants $I(R,E)$ which are
homogeneous of degrees $d_R$ and $d_E$ in $R$ and $E$, respectively.  If
$d_E=0$, we may regard $I(R,E)$ as an invariant $I(R)$ of $\calH_1$.  For
$I(R)$, we can apply Theorem C of \cite{BEG} and express it as a linear
combination of elementary invariants of $\calH_1$.  We thus assume
$d_E\ge1$, and again follow the arguments of \cite{BEG}.

The first step of the proof is to express $I(R,E)$ as a component of a
linear combination of partial contractions.  We denote by $\odot^{p,q}W^*$
the space of bisymmetric tensors of type $(p,q)$ on $W^*$ and by
$\odot_0^{p,q}W^*$ the subspace of $\odot^{p,q}W^*$ consisting of
trace-free tensors.  Let $e^*$ be the row vector 
$(0,\dots,0,1)\break\in W^*\otimes\sigma_{1,0}$. Then we have: \enddemo

\proclaim{Lemma}\label{sm}  
For some integer $m\le w-d_R-(n+1)d_E${\rm ,} a map
$$ 
  C\colon\calH_1\oplus\calH (-n-1)\to\odot{}_0^{m,m}W^*
  \otimes\sigma_{m-w,m-w}
$$  
is defined by making a linear combination of partial contractions of the
tensors $R^{(p,q)}, E^{(p,q)}$ and $e^*${\rm ,} $\c e^*$ such that
\begin{equation}
  C_{n\cdots n\,\c n\cdots\c n}=I.
  \label{CI}
\end{equation}
\endproclaim 

\demo{Proof}  
Since $R_{\alpha\c\beta}$ and $E_{\alpha\c\beta}$ satisfy 
\begin{equation}\label{link0}
  \begin{array}{ll}
    R_{\alpha0\c\beta}=(1-|\alpha|)R_{\alpha\c\beta},&\
    R_{\alpha\,\c{\beta}\,\c0}=(1-|\beta|) R_{\alpha\,\c{\beta}},
    \\
    E_{\alpha0\c\beta}=(-n-1-|\alpha|)E_{\alpha\c\beta},&\  
    E_{\alpha\,\c{\beta}\,\c0}=(-n-1-|\beta|) E_{\alpha\,\c{\beta}},
  \end{array} \hskip.5in
\end{equation} 
we can write $I(R,E)$ as a polynomial in the components of the form
$$
  {\wh R}^{k\c l}_{\alpha\c\beta}=
    R_{\alpha \underbrace{\scriptstyle{ n\cdots n}}_{k}\c\beta\,
  \underbrace{\scriptstyle{ \c n\cdots \c n}}_{l}},\quad
  {\wh E}^{k\c l}_{\alpha\c\beta}=
    E_{\alpha \underbrace{\scriptstyle{ n\cdots n}}_{k}\c\beta\,
  \underbrace{\scriptstyle{ \c n\cdots \c n}}_{l}},
$$ 
where $\alpha,\beta$ are lists of indices between 1 and $n-1$. We now regard
${\wh R}^{k\c l}_{\alpha\c\beta},{\wh E}^{k\c l}_{\alpha\c\beta}$ as tensors
on $\C^{n-1}$ by fixing $k,\c l$. For these tensors,  the Levi factor
$$ 
  L=
  \left\{
    \left(
      \begin{array}{ccc}
        \lambda & 0 & 0 \\
        0       & u & 0 \\
        0       & 0 & 1/\c\lambda
       \end{array}
     \right):
     u\in \U(n-1),\ \lambda\c\lambda^{-1}\det u=1
  \right\}
$$  
of $H$ acts as the usual tensorial action of $u$ up to a scale factor depending
on $\lambda$. Thus we may regard $I(R,E)$ as a $\U(n-1)$-invariant
polynomial.  By Weyl's classical invariant theory for $\U(n-1)$, we then see
that $I$ is expressed as a linear  combination of complete contractions of
the tensors 
${\wh R}^{k\c l}_{\alpha\c\beta},{\wh E}^{k\c l}_{\alpha\c\beta}$  for the
standard metric $\delta^{i\,\cj}$ on $\C^{n-1}$.  (See Lemma 7.4 of
\cite{BEG} for details of this discussion.)

We next replace the contractions with the metric $\delta^{i\,\cj}$ by 
those with the metric $g_0$. This is done by using the relation
$$
  \sum_{j=1}^{n-1}T_{j\,\cj}=-
  \sum_{i,j=0}^{n}g_0^{i\,\cj}T_{i\,\cj}+T_{0\,\c n}+T_{n\c 0}
  \qtext{for} 
  (T_{i\,\cj})\in W^*\otimes\c W^*.
$$  
Then several $0$ and $\c0$ come out as indices.  These can be eliminated by
using \eqref{link0}, and there remain only $n$ and $\c n$ as indices.  We
then get an expression of $I$ as a linear combination
$$
  I=\sum_{j=1}^k c_j C^{(j)}_{
  \underbrace{\scriptstyle{ n\cdots n}}_{p_j}
  \underbrace{\scriptstyle{\c n\cdots\c n}}_{q_j}},
$$ 
where each $C^{(j)}\in\bT^{p_j,q_j}_w$ is given by partial contraction of
the tensor products of several $R^{(p,q)}$ and $E^{(p,q)}$. In general,
$p_j,q_j$  $(1\le j\le k)$ are different.  Denoting by $m$ the maximum of
$p_j,q_j$  $(1\le j\le k)$, we define a tensor
$$
  C'=\sum_{j=1}^k c_j (\otimes^{m-p_j}e^*)\otimes
  C^{(j)}\otimes(\otimes^{m-q_j}\c e^*)\in \bT^{m,m}_w.
$$ 
Then we have $I=C'_{n\cdots n\c n\cdots\c n}$ because $e^*_n=\c e^*_n=1$. 
The map $C$ is now given by taking the trace-free bisymmetric part of $C'$.

To obtain the estimate $m\le w-d_R-(n+1)d_E$, we note that $C$ contains at
least one partial contraction which has no $e^*$ or no $\c e^*$. If such a
term consists of $R^{(p_j,q_j)}$, $E^{(p'_j,q'_j)}$ and several $e^*$
(resp.~$\c e^*$),  then $C$ takes values in
$\odot^{m,m}_0 W^*\otimes\sigma_{\kappa,\kappa}$ with
$\kappa=\sum_{j=1}^{d_R}(1-q_j)+\sum_{j=1}^{d_E}(-n-1-q'_j)$ (resp.~the
same relation with $p$ in place of $q$).  Hence noting $\kappa=m-w$ and
$p_j,q_j\ge2$, we obtain
$m-w\le\sum_{j=1}^{d_R}(-1)+\sum_{j=1}^{d_E}(-n-1)$, which is equivalent to
the desired estimate for $m$.
\enddemo

Now we regard $\calH_1\oplus\calH (-n-1)$ as the space of pairs of formal
power series $(\varphi,\eta)$ about $e_0\in W$ and write $I(\varphi,\eta)$ 
for $I(R,E)$ and $C(\varphi,\eta)$ for $C(R,E)$. If we replace the tensors
$R_{\alpha\c\beta}$ (resp.~$E_{\alpha\c\beta}$) in the partial contractions
in $C$ by the formal power series
$\pa_\zeta^\alpha\pa_{\c\zeta}^\beta\varphi$
(resp.~$\pa_\zeta^\alpha\pa_{\c\zeta}^\beta\eta$), we obtain a formal power
series about $e_0\in W$ which takes values in $\odot{}^{m,m}_0W^*$.
Restricting this power series to $\calQ$ and raising all indices by using
$g_0$, we obtain a map
\begin{eqnarray*}
&&  \wt C\colon\calH_1\oplus\calH(-n-1)
  \to
  \odot{}^{m,m}_0W\otimes\calJ(m-w)
\\
\noalign{\hbox{which satisfies $\wt C(\varphi,\eta)|_{e_0}=C(\varphi,\eta)$ when all
indices are raised.}}
\end{eqnarray*}

\vglue-20truept
Note that $\calH (k)$, $\calH^k$, $\calH_k$ and $\calJ(k)$ admit a natural
structure of $(\su(g_0),H)$-modules, where $\su(g_0)$ is the Lie algebra of
$\SU(g_0)$. For $\calH (k)$, $\calH^k$ and $\calJ(k)$ there are natural
$(\su(g_0),H)$-actions induced from the action of $\SU(g_0)$ on $W$ and
$\calQ$. For $\calH_k$,  a $(\su(g_0),H)$-action is induced via the
$H$-isomorphism $\calH_k\cong\calH (k)/\calH^k$.  We also consider the
complexification of these spaces and denote them, e.g., by $\calH^\C (k)$,
$\calJ^\C(k)$. Now we have:

\proclaim{Lemma}  {\rm (i)}  
There exists a unique $(\su(g_0),H)${\rm -}\/equivariant map 
$$
  \wt I\colon{\calH_1}\oplus\calH(-n-1)\to\calJ^\C(-w)
$$
such that $\wt I(\varphi,\eta)|_{e_0}=I(\varphi,\eta)$ for any 
$(\varphi,\eta)\in{\calH}_1\oplus\calH (-n-1)${\rm .}

{\rm (ii)}  
For any $(\varphi,\eta)\in\calH_1\oplus\calH (-n-1)${\rm ,}
\begin{equation}
  \wt C(\varphi,\eta)^{\alpha\c\beta}=
      \zeta ^{\alpha_1}\cdots  \zeta ^{\alpha_m}
    \c\zeta{}^{\beta_1}\cdots\c\zeta{}^{\beta_m}
  \wt I(\varphi,\eta).
  \label{wtCI}
\end{equation}
\endproclaim 

The proof of this lemma goes exactly the same way as those of Propositions
8.1 and 8.5 of \cite{BEG}, where the $(\su(g_0),H)$-equivariance of 
$\wt C$ is used essentially.

The final step in the proof of Theorem $5''$ is to obtain an explicit
expression of $\wt I$ in terms of $\wt C$  by differentiating both sides of
the equation \eqref{wtCI}.  We first introduce  differential operators 
$D_{i\,\cj}\colon\calE^\C (s+1)\to\calE^\C (s)$ for $(n+2s)(n+2s+1)\break\ne~0$ by
$$ 
  D_{i\,\cj} f=
  \left(\pa_{\zeta^i}-\frac{\zeta_i\Delta}{n+2s}\right)
  \left(\pa_{{\c \zeta}^j}-\frac{\c\zeta_j\Delta}{n+2s+1}\right)f,
$$   
where the index for $\zeta$ is lowered with $g_0$.  Then one can easily
check the following facts: (i) $D_{i\,\cj}(\mu\,f)=O(\mu)$, so that
$(D_{i\,\cj}f)|_{\cal Q}$ depends only on $f|_{\cal Q}$; (ii) For any
$f\in\calE^\C(s)$, 
$$ 
  D_{i\,\cj}(\zeta^i\c\zeta{}^j f)=c_s\,f,
  \qtext{where} 
  c_s=\frac{(n+s)^2(n+2s+2)}{n+2s}
$$  
and the repeated indices are summed over $0,1,\dots,n$; see Lemma 8.7
of \cite{BEG}. In view of (i), (ii) and taking an arbitrary extension of
$\wt I(\varphi,\eta)$ to $\calE (-w)$, we get 
$$
  \begin{array}{rl} 
    D_{\alpha_1\c \beta_1}D_{\alpha_2\c \beta_2}
    &\cdots
      D_{\alpha_m\c\beta_m}\wt C^{\alpha\c\beta}(\varphi,\eta)\\
    &=D_{\alpha_1\c\beta_1}D_{\alpha_2\c \beta_2}\cdots
      D_{\alpha_m\c\beta_m}
      (\zeta^{\alpha}\c\zeta{}^{\beta}\wt I(\varphi,\eta))\\ 
    &=c_{-w}c_{-w+1}\cdots c_{-w+m-1}\,\wt I(\varphi,\eta)
      \qtext{on}\calQ.
  \end{array}
$$  
Since $m-w\le-d_R-(n+1)d_E$ and $1\le d_E$, we have $m-w\le-n-1$. Hence,
all $D_{i\cj}$ appearing above are well-defined  and all $c_s\ne 0$.
Therefore,
$$ 
  I(\varphi,\eta)=
  \frac{1}{c_{-w}\cdots c_{-w+m-1}} 
  D_{\alpha_1\c \beta_1}\cdots D_{\alpha_m\c\beta_m}
  \wt C^{\alpha\c\beta}(\varphi,\eta)|_{e_0},
$$  
and $I$ is expressed as a linear combination of complete contractions.

\numbereddemo{{R}emark} 
The assumption $d_E\ge1$ was used only in the final step of the proof to
ensure $c_s\ne0$. The argument above is valid even if $d_E=0$ as long as
$d_R\ge n$. This is exactly the proof of Theorem C of \cite{BEG} for
invariants of high degrees. To treat the invariants of low degrees on
$\calH_1$, the authors used an entirely different argument.
\enddemo

\numbereddemo{{R}emark} \label{rem-graham} The tensors $E^{(p,q)}$ used in this
section are modeled on the biholomorphically invariant tensors
$$
  E^{(p,q)}_k=\nabla{}^p\c\nabla{}^q(|z^0|^{-2k(n+1)}\eta_k),
$$ 
which were introduced by Graham \cite{G2}. He used these tensors  to
construct CR invariants from the complete contractions of the form
$$
  \contr\big(R^{(p_1,q_1)}\otimes\cdots R^{(p_d,q_d)}\otimes
  E_{k_1}^{(p_1',q_1')}\otimes\cdots E_{k_{d'}}^{(p'_{d'},q'_{d'})}\big).
$$ 
Such complete contractions give rise to CR invariants if, for example,
$p_j,q_j< n+2$ and $p'_j,q_j'<n+1$.  This class of CR invariants correspond
to $\calC$-independent Weyl invariants which contain the covariant
derivatives of the Ricci tensor.  In fact, $g[r]$ is Ricci flat if and only
if $\eta_k=0$ for all $k\ge1$, because the Ricci form of $g[r]$ is given by
$\pa\c\pa\log J[r]$. Thus, a CR invariant depending on $E^{(p,q)}_k$  must
contain the covariant derivatives of the Ricci tensor when  it is expressed
as a Weyl invariant. 
\enddemo

\section{Proof of Theorem 3}

By virtue of Theorem 2, it suffices to prove:

\proclaim{Proposition}\label{prop-c-dependence} 
Let $n\ge3$ $($resp. $n=2)${\rm .} Then every Weyl invariant of weight $w\le n+2$
$($resp. $w\le5)$ is $\calC$\/{\rm -}\/independent{\rm .}  For $w=n+3$ $($resp. $w=6)${\rm ,}
there exists a $\calC$\/{\rm -}\/dependent Weyl invariant of weight $w${\rm .}
\endproclaim

{\it Proof of Proposition $\ref{prop-c-dependence}$}. 
Take a Weyl polynomial $W_\#$ of weight $w$ and set $I(A,C)=I_W(A,C)$. We
begin by inspecting the linear part of $I(A,C)$.

\proclaim{Lemma}\label{lem-linear-term} 
If $I(A,C)$ has nonzero linear part{\rm ,} then $w=n+2$ and the linear part is a
constant multiple of $\Delta_x^{n+2}f_A(e_0)$ with $f_A$ as in {\rm \eqref{bcA}.}
\endproclaim 

\demo{Proof} If $W_\#(R)$ has no linear terms, neither does $I(A,C)$. Thus it
suffices to consider the case where  $W_\#(R)$ is a linear complete contraction
$\contr(R^{(p,p)})$. In this case, \eqref{R-linear} implies
\begin{equation}\label{linear-weyl}
  I(A,C)=\Delta^p\varphi(e_0)+Q(A,C),
\end{equation}  
where $Q(A,C)$ is a polynomial in
$(A_{\alpha\c\beta}^l,C_{\alpha\c\beta}^l)$ without linear terms.  By
\eqref{eq-varphi} and \eqref{Delta-Delta}, $\Delta^p\varphi(e_0)=0$ if
$p\ne n+2$ and $\Delta^{n+2}\varphi(e_0)=-\Delta_x^{n+2}f_A(e_0)$. Thus we
obtain the lemma.
\enddemo

We next consider nonlinear terms in $I(A,C)$. Since $I(A,C)$ is invariant
under the action of $\U(n-1)$, it is written as a linear combination of
complete contractions of the form
\begin{equation}\label{contrAC}
  \contr'
  \left(
    \bA^{l_1}_{p_1\c q_1}\otimes\cdots\otimes
    \bA^{l_d}_{p_d\c q_d}\otimes {\bf C}^{l'_1}_{p'_1\c q'_1}
    \otimes\cdots\otimes{\bf C}^{l'_{d'}}_{p'_{d'}\c q'_{d'}}
  \right)
\end{equation}  
with 
$\sum_{j=1}^d(p_j+q_j+2l_j-2)+\sum_{j=1}^{d'}(p'_j+q'_j+2l'_j+2n+2)=2w$.
Here $\bC^l_{p\c q}=(C_{\alpha\c\beta}^l)_{|\alpha|=p,|\beta|=q}$ is
regarded as a tensor of type $(p,q)$ on $\C^{n-1}$ and the contraction is
taken with respect to $\delta^{i\cj}$ for some pairing of lower indices.
Suppose \eqref{contrAC} is nonlinear and contains the variables
$C_{\alpha\c\beta}^l$, so that $d+d'\ge2$ and $d'\ge1$. Then $p_j+q_j\ge4$
implies $w\ge n+2$.  The equality $w=n+2$ holds only for
$\contr'(\bA^0_{2\c2}\otimes{\bC}^0_{0\c0})$, while by (N2),
$$
  \contr'(\bA^0_{2\c2}\otimes{\bC}^0_{0\c0}) =\contr'(\bA^0_{2\c2})C^0=0,
$$ 
where $C^0$ is the only one component of $\bC^0_{0\c0}$. Thus $I(A,C)$
containing $C_{\alpha\c\beta}^l$ has weight $\ge n+3$. In case $w=n+3$,
there are only two types of contractions of the form 
\eqref{contrAC}, namely,
\begin{equation}\label{contrAAC}
  \contr'(
   \bA^0_{2\c 2}\otimes\bA^0_{2\c 2}
   \otimes\bC^0_{0\c0})\quad\hbox{and}\quad
  \contr'(\bA^l_{p\c q}\otimes\bC^{l'}_{p'\c q'}),
\end{equation} 
where $p+p'=q+q'=3-l-l'$.  The contractions of the second type always vanish
by (N2) (see \S 3), and the first ones vanish except for the case 
$\|\bA^0_{2\c2}\|^2\,C^0 =
\sum_{i,j,k,l=1}^{n-1}|A_{ij\,\c k\,\c l}^0|^2\,C^0$; 
this also vanishes for $n=2$ because 
$A_{11\,\c1\c1}^0=\tr\,\bA^0_{2\c2}=0$. This completes the proof of the
first statement of Proposition \ref{prop-c-dependence}.

To prove the second statement, we consider, for $n=2$, a complete
contraction of weight $6$:
$$
  W_{2}=\sum_{|\alpha|=6,|\beta|=2} R_{\alpha\c\beta}R^{\c\beta\alpha}
$$ 
and, for $n\ge 3$, 
$$
  W_{n}=\sum_{|\alpha|=|\beta|=2,|\gamma|=n+2}
   R_{\alpha\,\c\beta}R^{\c\beta}{}^\gamma R_\gamma{}^\alpha,
$$ 
which has weight $n+3$. Here indices are raised by using $g_0$. These
complete contractions give $\calC$-dependent Weyl invariants.  In fact:

\proclaim{Lemma}
Let $I_n(A,C)=I_{W_n}(A,C)${\rm .} Then
\begin{equation}\label{p2}
  I_2(A,C)=72\cdot6!\,(C^0)^2+Q_2(A,C),
\end{equation} 
where $Q_2$ is a polynomial in $(A_{\alpha\c\beta}^l,C_{\alpha\c\beta}^l)$
such that $Q_2(0,C)=0$. For $n\ge3${\rm ,}
\begin{equation}\label{pn}
  I_n(A,C)=(-1)^n 64\,(n+2)!\,\|\bA_{2\c 2}^0\|^2 C^0+Q_n(A),
\end{equation} 
where $Q_n(A)$ is a polynomial in $A_{\alpha\c\beta}^l${\rm .}
\endproclaim 

\demo{Proof}  We first prove \eqref{pn}.  Since
$I_n(A,C)=c_n\,\|\bA_{2\c2}^0\|^2 C^0+Q_n(A)$ for some constant $c_n$, to
determine $c_n$, we consider a family of  surfaces with real parameter $s$
$$
  2u=|z'|^2+f_s(z',\c z'),\qtext{where}
  f_s(z',\c z')=2s\,\Re\,(z^1)^2(\c z^2)^2.
$$ 
Let $A_s\in\calN$ denote the list of normal form coefficients of this
surface and $C_t\in\calC$ the element such that $C^0=t$ and  all the other
components vanish. Then
\begin{equation}\label{Pn-1}
  I_n(A_s,C_t)=c_n\,s^2t+Q_n(A_s)=c_n\,s^2t+O(s^{n+3}).
\end{equation} 
On the other hand, we see by \eqref{R-linear} that the components
$R_{\alpha\c\beta}(s,t)$ of $\iota(A_s,C_t)$ satisfy
$$
  R_{\alpha\c\beta}(s,t)=S_{\alpha\c\beta}+O(s^2+t^2),
$$ 
where $S_{\alpha\c\beta}=\pa_\zeta^\alpha\pa_{\c\zeta}^\beta\varphi(e_0)$
with $\varphi=\calL(f_s,t)$. Thus 
\begin{equation}\label{Pn-2}
  I_n(A_s,C_t)=W'_n+O((s^2+t^2)^2),
\end{equation}  
where
\begin{equation}\label{Wn} 
  W'_n=\sum_{|\alpha|=|\beta|=2,|\gamma|=n+2}
  S_{\alpha\c\beta}\,S^{\c\beta\gamma}\,S_\gamma{}^\alpha.
\end{equation} 
Comparing \eqref{Pn-1} with \eqref{Pn-2}, we get
$$
  W_n'=c_n\,s^2 t.
$$ 
Since 
$$
  \varphi=-|\zeta^0|^2 f_s+t\,|\zeta^0|^{-2(n+1)}\mu^{n+2}/(n+2)!,
$$  
the term $S_{\alpha\c\beta}$ in the sum of \eqref{Wn}  vanishes except for
$S_{11\c2\c2}=S_{22\c1\c1}=-4s$. Using
$$
  \sum_{|\gamma|=n+2}
  S^{\c2\c2\,\gamma}S_\gamma{}^{11}=
  \sum_{|\gamma|=n+2}
  S_{\c1\c1\,\gamma}S^\gamma{}_{22}=
  \c{\sum_{|\gamma|=n+2}
  S^{\c1\c1\,\gamma}S_\gamma{}^{22}},
$$ 
we then obtain
\begin{eqnarray*}
  W'_n&=&-4s\sum_{|\gamma|=n+2}
  \left(
    S^{\c2\c2\,\gamma}S_\gamma{}^{11}+
    S^{\c1\c1\,\gamma}S_\gamma{}^{22}
  \right)\\
  &=&-8s\,\Re\sum_{|\gamma|=n+2}
    S_{\gamma\,\c1\c1}S^\gamma{}_{22}.
\end{eqnarray*} 
In the last sum, $S_{\gamma\,\c1\c1}$ vanishes unless $\gamma$ is a
permutation of $0\cdots022$ or $11n\cdots n$, while
$$
  S_{0\cdots022\,\c1\c1}=
  S_{22}{}^{11n\cdots n}=(-1)^{n+1}4\cdot n!\,s,
$$
$$
  S_{11n\cdots n\,\c1\c1}=S_{22}{}^{0\cdots022}=2\,t.
$$ 
Therefore $W_n'=(-1)^n 64\, (n+2)!\,s^2\,t$, so that
$c_n=(-1)^n 64\,(n+2)!$.

We next prove \eqref{p2}. Since $I_2(A,C)$ contains no linear term, we have
$$
  I_2(A,C)=c_2\,(C^0)^2+Q_2(A,C)
$$ 
for a constant $c_2$. Restriction of this formula to $(A,C)=(0,C_t)$ yields
\begin{equation}\label{P2t}
  I_2(0,C_t)=c_2\,t^2+O(t^3),
\end{equation} 
while, by the expression $\varphi=\calL(0,t)=t\,|\zeta^0|^{-6}\mu^4/4!$,
$$
  I_2(0,C_t)=W_2'+O(t^3),\quad
  W'_2=\sum_{|\alpha|=6,|\beta|=2}
  S_{\alpha\,\c\beta}\,S^{\c\beta\,\alpha}.
$$ 
Since $S_{\alpha\c0\c k}=S^{\c2\,\c k\,\alpha}=0$ for $k=0,1,2$ and any
list $\alpha$, we have
$$
  W'_2=\sum_{|\alpha|=6}S_{\alpha\,\c1\c1}\,S^{\c1\c1\,\alpha}.
$$ 
In this sum, $S_{\alpha\,\c1\c1}$ vanishes unless $\alpha$ is a permutation
of $001122$, while 
$$
  S_{001122\c1\c1}=S^{\c1\c1\,001122}=4!\,t.
$$  
Thus $W'_2=72\cdot 6!\,t^2$. This together with \eqref{P2t} yields
$c_2=72\cdot 6!$.
\enddemo

\numbereddemo{{R}emark}\rm 
As a consequence of Lemma \ref{lem-linear-term}, we see that a CR invariant
$I(A)$ of weight $w$ can contain linear terms only when $w=n+1$ and that the
linear part must be a constant multiple of $\Delta_x^{n+2}f_A(e_0)$. This
fact is equivalent to Theorem 2.3 of Graham \cite{G2}. 
\enddemo

\AuthorRefNames  [HKN2]

\end{document}